\documentclass[12pt]{article}
\usepackage{amssymb,a4}
\usepackage{amsmath,amsfonts,amssymb,amsthm}
\input amssym.def
\input amssym.tex
\def \wt{{\rm wt}}

\def \End{{\rm End}}
\def \Hom{{\rm Hom}}

\def\Res{{\rm Res}}
\def\wt{{\rm wt}}

\def\de{\delta}

\def\C{{\mathbb C}}

\def\Z{{\mathbb Z}}

\def \Q{\mathbb Q}
\def\1{{\bf 1}}
\def\l{{\lambda}}
\def \pf{\noindent {\bf Proof: \,}}
\def\theequation{5.\arabic{equation}}

\renewcommand{\theequation}{\thesection.\arabic{equation}}

\newtheorem{theorem}{Theorem}[section]
\newtheorem{prop}[theorem]{Proposition}
\newtheorem{lem}[theorem]{Lemma}
\newtheorem{coro}[theorem]{Corollary}
\newtheorem{remark}[theorem]{Remark}
\theoremstyle{definition}
\newtheorem{definition}[theorem]{Definition}

\begin{document}
\begin{center}
{\Large {\bf Bimodules and $g$-rationality of vertex operator algebras}} \\

\vspace{0.5cm} Chongying Dong\footnote{Supported by NSF grants,
China NSF grant 10328102 and  a Faculty research grant from  the
University of California at Santa Cruz.}
\\
Department of Mathematics\\ University of
California\\ Santa Cruz, CA 95064 \\
Cuipo Jiang\footnote{Supported  by China NSF grant
10571119.}\\
 Department of Mathematics\\ Shanghai Jiaotong University\\
Shanghai 200030 China\\
\end{center}
\hspace{1cm}
\begin{abstract}  This paper studies the twisted representations of vertex
operator algebras. Let $V$ be a vertex operator algebra and $g$ an
automorphism of $V$ of finite order $T.$ For any
$m,n\in\frac{1}{T}\Z_+$, an $A_{g,n}(V)$-$A_{g,m}(V)$-bimodule
$A_{g,n,m}(V)$ is constructed. The collection of these bimodules
determines any admissible $g$-twisted $V$-module completely. A
Verma type admissible $g$-twisted $V$-module is constructed
naturally from any $A_{g,m}(V)$-module. Furthermore, it is shown
with the help of bimodule theory that a simple vertex operator
algebra $V$ is $g$-rational if and only if its twisted associative
algebra $A_g(V)$ is semisimple and each irreducible admissible
$g$-twisted $V$-module is ordinary.
\end{abstract}

 \section{Introduction}
\def\theequation{1. \arabic{equation}}
\setcounter{equation}{0}

This paper deals with twisted representations of vertex operator algebras
using the ideas of bimodules developed in \cite{DJ1}-\cite{DJ3}. The main
result is a characterization of twisted rationality in terms of semisimplicity
of cerain associative algebra defined and studied in \cite{DLM2}.

Twisted representations which are also called twisted sectors or
twisted modules are the main ingredients in orbifold conformal
field theory (see \cite{DHVW1}-\cite{DHVW2}, \cite{L1}-\cite{L2},
\cite{FLM1}-\cite{FLM2}, \cite{DVVV}, \cite{DM}, \cite{DLM0},
\cite{HMT}, \cite{DLM2}, \cite{DLM4}, \cite{DLM5}, \cite{DY} and \cite{MT}).
The twisted sectors  play a fundamental role in the construction
of the moonshine vertex operator algebra $V^{\natural}$
\cite{FLM2} and other orbifold vertex operator algebras
\cite{DGM}. Although there is a lot of progress in the study of
twisted sectors and orbifold conformal field theory, the
semisimplicity of various twisted module categories has not been
understood fully.

Let $V$ be a vertex operator algebra and $g$ an automorphism of
finite order $T.$  There are three different notions of
$g$-twisted modules. That is,  weak $g$-twisted modules,
admissible $g$-twisted modules and  ordinary $g$-twisted modules
(see \cite{FFR}, \cite{D}, \cite{DLM2}). An ordinary $g$-twisted
$V$-module is admissible and an admissible $g$-twisted $V$-module
is a weak $g$-twisted $V$-module. They differ by some grading
assumptions. The main axiom in these modules is the twisted Jacobi
identity  which was motivated
 by the twisted vertex operator operators studied in
\cite{L1}-\cite{L2} and \cite{FLM1}-\cite{FLM2}.

We call a vertex operator algebra $V$ $g$-rational if the
admissible $g$-twisted $V$-module category is semisimple. It is
proved in \cite{DLM2} that if $V$ is $g$-rational then there are
only finitely many irreducible admissible $g$-twisted $V$-modules
up to isomorphism and each irreducible admissible $g$-twisted
$V$-module is ordinary. So the concept of $g$-rationality is an
analogue of semisimplicity of associative algebras. In fact, the
$g$-rationality will be understood in terms of semisimplicity of
an associative algebra $A_g(V)$  investigated in \cite{DLM2}.

Stimulated by the $A(V)$-theory developed in \cite{Z}, an
associative algebra $A_g(V)$ is defined and studied in
\cite{DLM2}. In order to state the connection between the twisted
representation theory of $V$ and the representation theory of
$A_g(V)$ let $M=\bigoplus_{n\in\frac{1}{T}\Z_+}M(n)$ be an
admissible $g$-twisted $V$-module with $M(0)\ne 0.$ Then $M(0)$ is
an $A_g(V)$-module. Moreover, the map $M\to M(0)$ gives a one to
one correspondence between the irreducible admissible $g$-twisted
$V$-modules and simple $A_g(V)$-modules. These results reduce the
classification of irreducible admissible $g$-twisted $V$--modules
to the classification of simple $A_g(V)$-modules. So the
classification of irreducible admissible $g$-twisted $V$--modules
is settled down at least theoretically.

The main purpose of this paper is to establish a relationship
between the $g$-rationality of $V$ and the semisimplicity of
$A_g(V).$  It has already been proved in \cite{DLM2} that the
$g$-rationality of $V$ implies  the semisimplicity of $A_g(V).$ We
prove in this paper that $V$ is $g$-rational if and only if
$A_g(V)$ is semisimple and each irreducible admissible $g$-twisted
$V$-module is ordinary. Note that the $g$-rationality is an
external condition on $V.$ The new result essentially gives an
internal characterization of $g$-rationality as $A_g(V)$ is a
quotient of $V$ \cite{DLM2}. In the case that $g=1$ this result
has been obtained in \cite{DJ3}.

The main idea comes from \cite{DJ1}-\cite{DJ3}. The associative
algebra $A_{g}(V)$ was generalized to associative algebras
$A_{g,n}(V)$ for any $n\in\frac{1}{T}\Z_+$ so that
$A_{g,0}(V)=A_g(V)$ \cite{DLM4}. For an admissible $g$-twisted
$V$-module $M=\bigoplus_{n\in\frac{1}{T}\Z_+}M(n)$ with $M(0)\ne
0,$ $M(m)$ is a module for $A_{g,n}(V)$ for $m\leq n.$ So
$A_{g,n}(V)$ gives  more information on $M$ than $A_g(V).$ Most
importantly, $V$ is $g$-rational if and only if $A_{g,n}(V)$ is
semisimple for all $n\in\frac{1}{T}\Z_{+}$. Our approach is to
prove that if $A_{g}(V)$ is semisimple then $A_{g,n}(V)$ is
semisimple for all $n.$

We first construct $A_{g,n}(V)$-$A_{g,m}(V)$-bimodules
$A_{g,n,m}(V)$ (for $n,m\in\frac{1}{T}\Z_+$) which establish a
bridge between $A_{g}(V)$ and $A_{g,n}(V).$ These bimodules carry
the semisimplicity information from $A_{g}(V)$ to all
$A_{g,n}(V).$ From the point of view of representation theory, the
$A_{g,n}(V)$-$A_{g,m}(V)$-bimodule  $A_{g,n,m}(V)$ is a universal
covering of the  $A_{g,n}(V)$-$A_{g,m}(V)$-bimodule
$\Hom_{\C}(M(m),M(n))$ for any admissible $g$-twisted $V$-module
$M.$ The importance of the construction of these bimodules is that
it gives a concrete construction of the Verma type admissible
$g$-twisted $V$-module $M(U)$ generated by an $A_{g,m}(V)$-module
$U$ such that $M(U)(n)=A_{g,n,m}(V)\otimes_{A_{g,m}(V)}U$ for all
$n\in\frac{1}{T}\Z_+.$ Using this construction we can prove that
there is a natural invariant pairing between $M(U^*)$ and $M(U)$
such that the right radical of this pairing is exactly the maximal
proper submodule of $M(U)$ if $U$ is irreducible. The $M(U^*)$ is
an admissible $g^{-1}$-twisted $V$-module as $U^*$ is an
$A_{g^{-1},m}(V)$-module instead of $A_{g,m}(V)$-module. So this
makes the admissible $g$-twisted module theory more comparable
with the classical highest weight module theory for affine Lie
algebras or the Virasoro algebra. With the help of this
construction we can also prove that if $A_g(V)$ is semisimple,
then the Verma type admissible $g$-twisted $V$-module $M(U)$
generated by an irreducible $A_g(V)$-module $U$ is irreducible.
This is the key step in the proof of the main theorem.

Since the setting and most results in this paper are modelled on
those in \cite{DJ1} and \cite{DJ3} which deals with the case
$g=1,$ we omit a lot of details in this paper and refer the reader
to \cite{DJ1} and \cite{DJ3}.

 \section{The associative algebra $A_{g,n}(V)$}
\def\theequation{2.\arabic{equation}}
\setcounter{equation}{0}
Let $(V,Y,{\1},\omega)$ denote, as
usual, a vertex operator algebra as defined in [FLM2] (see also [B])
and $g$ be an automorphism of $V$ of finite order
$T.$  Decompose $V$ into a direct sum of eigenspaces of $g$
\begin{equation}\label{g2.1}
V=\bigoplus_{r\in \Z/T\Z}V^r
\end{equation}
where $V^r=\{v\in V|gv=e^{-2\pi ir/T}v\}$. We first review the
weak, admissible and ordinary $g$-twisted modules from \cite{DLM2}
(see also \cite{FLM2}, \cite{FFR} and \cite{D}).

\begin{definition} A {\em weak $g$-twisted $V$-module} $M$ is a vector space equipped
with a linear map
\begin{eqnarray*}
Y_{M}(\cdot,z): & &V\to (\End\,M)\{z\}\\
& &v\mapsto\displaystyle{ Y_M(v,z)=\sum_{n\in\Q}v_nz^{-n-1}\ \ \
(v_n\in \End\,M)}
\end{eqnarray*}
which satisfies the following conditions for all $0\leq r\leq
T-1,$ $u\in V^r$, $v\in V,$ $w\in M$,
\begin{eqnarray*}
& &Y_M(u,z)=\sum_{n\in \frac{r}{T}+\Z}u_nz^{-n-1} \label{1/2}\\
& &u_lw=0\ \ \
\mbox{for}\ \ \ l>>0\label{vlw0}\\
& &Y_M({\1},z)={\rm id}_{M};\label{vacuum}
\end{eqnarray*}
\begin{eqnarray*}
\displaystyle{z^{-1}_0\delta\left(\frac{z_1-z_2}{z_0}\right)
Y_M(u,z_1)Y_M(v,z_2)-z^{-1}_0\delta\left(\frac{z_2-z_1}{-z_0}\right)
Y_M(v,z_2)Y_M(u,z_1)}\\
\displaystyle{=z_2^{-1}\left(\frac{z_1-z_0}{z_2}\right)^{-r/T}
\delta\left(\frac{z_1-z_0}{z_2}\right) Y_M(Y(u,z_0)v,z_2)}.
\end{eqnarray*}
\end{definition}

As mentioned in \cite{DLM2} (see also \cite{FLM2}), the twisted
Jacobi identity is equivalent to the associativity formula
\begin{eqnarray}\label{ea}
(z_{0}+z_{2})^{k+\frac{r}{T}}Y_{M}(u,z_{0}+z_{2})Y_{M}(v,z_{2})w
=(z_{2}+z_{0})^{k+\frac{r}{T}}Y_M(Y(u,z_0)v,z_2)w
\end{eqnarray}
where $w\in M$ and $k$ is a nonnegative integer such that
$z^{k+\frac{r}{T}}Y_{M}(u,z)w$ involves only nonnegative integral
powers of $z,$ and the commutator formula
\begin{eqnarray}\label{cea}
& &\ \ \ \  [Y_{M}(u,z_{1}),Y_{M}(v,z_{2})]\nonumber\\
& &=\Res_{z_{0}}z_2^{-1}\left(\frac{z_1-z_0}{z_2}\right)^{-r/T}
\delta\left(\frac{z_1-z_0}{z_2}\right)Y_M(Y(u,z_0)v,z_2).\label{ec}
\end{eqnarray}


\begin{definition}\label{d3.2} An ordinary $g$-{\em twisted $V$-module} is
a weak $g$-twisted $V$-module $M$ which carries a $\C$-grading
induced by the spectrum of $L(0).$ That is,
$$
M=\bigoplus_{\lambda \in{\C}}M_{\lambda}
$$
where $M_{\l}=\{w\in M|L(0)w=\l w\},$ where $L(0)$ is a component
operator of $Y_M(\omega,z)=\sum_{n\in\Z}L(n)z^{-n-2}.$
 Moreover we require that
$ M_{\l}$ is finite dimensional and for fixed $\l,$ $M_{\frac{n}{T}+\l}=0$
for all small enough integers $n.$
\end{definition}

Let $\Z_+$ be the set of nonnegative integers.
\begin{definition} An {\em admissible} $g$-twisted $V$-module
is a  weak $g$-twisted $V$-module $M$ which carries a $\frac{1}{T}{\Z}_{+}$-grading
$$
M=\bigoplus_{n\in\frac{1}{T}\Z_+}M(n)
$$
which satisfies the following
$$
v_{m}M(n)\subseteq M(n+\wt v-m-1)
$$
for homogeneous $v\in V.$ \end{definition}

It is easy to show that an ordinary $g$-twisted $V$-module is
admissible. If $g=1$ we get the weak, ordinary and admissible
$V$-modules.

We say that $V$ is $g$-{\em rational} if every admissible
$g$-twisted $V$-module is completely reducible. $V$ is called {\em
rational} if $V$ is $1$-rational. It is proved in \cite{DLM2} that
if $V$ is $g$-rational then there are only finitely many
irreducible admissible $g$-twisted $V$-modules up to isomorphism
and each irreducible admissible module is ordinary.

Next we present the $A_{g,n}(V)$-theory following \cite{DLM4}.
Fix $n=l+\frac{i}{T}\in\frac{1}{T}\Z_{+}$ with $l$ a nonnegative
integer and $0\leq i\leq T-1.$ For $0\leq r\leq T-1$, define
$\delta_i(r)=1$ if $i\geq r$ and $\delta_i(r)=0$ if $i<r$.
 We also set $\de_{i}(T)=0$. Let
$O_{g,n}(V)$ be the linear span of all $u\circ_{g,n} v$ and
$L(-1)u+L(0)u$ where for homogeneous $u\in V^r$ and $v\in V,$ $$
u\circ_{g,n} v=\Res_{z}Y(u,z)v\frac{(1+z)^{\wt
u-1+\delta_i(r)+l+r/T}}{z^{2l +\delta_{i}(r)+\delta_{i}(T-r)+1 }}.
$$
We also define a second product $*_{g,n}$ on $V$ for $u\in V^r$
and $v$ as follows: $$ u*_{g,n}v=\sum_{m=0}^{l}(-1)^m{m+l\choose
l}\Res_zY(u,z)\frac{(1+z)^{\wt\,u+l}}{z^{l+m+1}}v
$$
if $r=0$ and $u*_{g,n}v=0$ if $r>0.$

Define the linear space $A_{g,n}(V)$ to be the quotient
$V/O_{g,n}(V).$ Then $A_{g,0}(V)=A_g(V)$ has been  defined and
studied in \cite{DLM2} already.

\begin{remark} The definition of
$u\circ_{g,n} v$ in \cite{DLM4} is not correct where it is defined
as
 $$
u\circ_{g,n} v=\Res_{z}Y(u,z)v\frac{(1+z)^{\wt
u-1+\delta_i(r)+l+r/T}}{z^{2l +\delta_{i}(r)+\delta_{i}(T-r)}}
$$
with $\delta_i(T)=1.$ But  the results and proofs in \cite{DLM4}
remain valid.
\end{remark}

Let $W$ be a weak $g$-twisted $V$-module and
$m\in\frac{1}{T}{\mathbb Z}_{+}.$ Following \cite{DLM4} we define
$$ \Omega_{m}(W)=\{w\in W|u_{{\wt}u-1+k}w=0, {\rm for \ all \
homogeneous } \ u\in V \ {\rm and} \ k>m\}.$$

The following theorem is obtained in \cite{DLM4}.

\begin{theorem}\label{t2.4}  Let $V$ be a vertex operator algebra and $g$
an automorphism of $V$ of finite order $T.$ Let
$M=\bigoplus_{m\in\frac{1}{T}\Z_+}M(m)$ be an admissible
$g$-twisted $V$-module. Let $n\in\frac{1}{T}\Z_+$. Then

 (1)
$A_{g,n}(V)$ is an associative algebra whose product is induced by
$*_{g,n}.$

(2) The identity map on $V$ induces an algebra epimorphism from
$A_{g,n}(V)$ to $A_{g,n\!-\frac{1}{T}}(V).$

(3) Let $W$ be a weak $g$-twisted $V$-module. Then $\Omega_n(W)$
is an $A_{g,n}(V)$-module such that $v+O_{g,n}(V)$ acts as
$o(v)=v_{\wt v-1}$ for homogeneous $v.$

(4) Each $M(m)$ for $m\leq n$ is an $A_{g,n}(V)$-submodule of
$\Omega_n(M).$ Furthermore, $M$ is irreducible if and only if each
$M(n)$ is an irreducible $A_{g,n}(V)$-module.

(5) For any $A_{g,n}(V)$-module $U$ which cannot factor through
$A_{g,n-\frac{1}{T}}(V)$ there is a unique Verma type admissible
$g$-twisted $V$-module $\bar{M}(U)$ generated by $U$ so that
$\bar{M}(U)(0)\ne 0$ and $\bar{M}(U)(n)=U.$ Moreover, for any weak
$g$-twisted $V$-module $W$ and any $A_{g,n}(V)$-module
homomorphism $f$ from $U$ to $\Omega_n(W)$ there is a unique
$V$-module homomorphism from $\bar{M}(U)$ to $W$ which extends
$f.$

(6) $V$ is $g$-rational if and only if $A_{g,n}(V)$ are finite
dimensional semisimple algebras for all $n\in\frac{1}{T}\Z_{+}.$

(7) If $V$ is $g$-rational then there are only finitely many
irreducible admissible $g$-twisted $V$-modules up to isomorphism
and each irreducible admissible $g$-twisted $V$-module is
ordinary.

(8) The linear map  $\phi$ from $V$ to $V$ defined by
$\phi(u)=e^{L(1)}(-1)^{L(0)}u$  for $u\in V$ induces
an anti-isomorphism from $A_{g,n}(V)$ to $A_{g^{-1},n}(V).$
\end{theorem}

\section{ $A_{g,n}(V)$-$A_{g,m}(V)$-bimodule $A_{g,n,m}(V)$}
\def\theequation{3.\arabic{equation}}
\setcounter{equation}{0}

Let $V=(V,Y,{\bf 1},\omega)$ be a vertex operator algebra, and let
$g$ be an automorphism of $V$ of finite order $T$. This section is
an extension of bimodule theory developed in \cite{DJ1} from the
untwisted case to the twisted case. In particular we will
construct an $A_{g,n}(V)$-$A_{g,m}(V)$-bimodule $A_{g,n,m}(V).$

For $k\in\Z$, we denote the image of $k$ in $\Z/T\Z$ by $\bar{k}$.
Without confusion, if $0\leq k\leq T-1$, we sometimes also denote
$\bar{k}\in\Z/T\Z$ by $k$. Let $m,p,n\in(1/T)\Z_{+}$. Then
$m=l_{1}+(i_{1}/T),$ $p=l_{2}+(i_{2}/T),$
$n=l_{3}+(i_{3}/T)\in(1/T)\Z$ with $l_{1}, l_{2}, l_{3}$ three
nonnegative integers and $0\leq i_{1},i_{2},i_{3}\leq T-1$. In the
following discussion,  we always denote $m,n,p$ as above until
further notice.

Recall the decomposition  (\ref{g2.1}).  For homogeneous $u\in
V^r$, $v\in V$, define product $*_{g,m,p}^{n}$ on $V$ as follows:
\begin{eqnarray*}
& &  u*_{g,m,p}^{n}v=\sum\limits_{i=0}^{l_{2}}
(-1)^{i}{l_{1}+l_{3}-l_{2}-1+\de_{i_{1}}(r)+\de_{i_{3}}(T-r)+i\choose
i}\\
& & \ \ \ \ \ \ \ \ \ \ \ \ \ \ \ \ \
\cdot\Res_{z}\frac{(1+z)^{{\wt}u-1+l_{1}+\delta_{i_{1}}(r)+r/T}}{z^{l_{1}+l_{3}-l_{2}+
\de_{i_{1}}(r)+\de_{i_{3}}(T-r)+i}}Y(u,z)v
\end{eqnarray*}
if   $\overline{i_{2}-i_{3}}=r$ and
$$u*_{g,m,p}^{n}v=0$$
otherwise.

 If $n=p$, we denote $\ast_{g,m,p}^{n}$ by
$\bar{\ast}_{g,m}^{n}$. In this case, $u\bar*_{g,m}^{n}v=0$ if
$r\neq 0$ and
$$
u\bar*_{g,m}^{n}v=\sum\limits_{i=0}^{l_{3}}
(-1)^{i}{l_{1}+i\choose
i}\Res_{z}\frac{(1+z)^{{\wt}u+l_{1}}}{z^{l_{1}+i+1}}Y(u,z)v
$$
for $r=0.$ One can easily check that ${\bf 1}\bar*_{g,m}^{n}u=u$, for $u\in
V^r$.

 If $m=p$, we denote $\ast_{g,m,p}^{n}$ by
$\ast_{g,m}^{n}$. In this case, $u\ast_{g,m}^{n}v=0$ if
$\overline{i_{1}-i_{3}}\neq  r$; if $\overline{i_{1}-i_{3}}=r$,
then $-1+\de_{i_{1}}(r)+\de_{i_{3}}(T-r)=0$. So
$$
u\ast_{g,m}^{n}v=\sum\limits_{i=0}^{l_{1}} (-1)^{i}{l_{3}+i\choose
i}\Res_{z}\frac{(1+z)^{{\wt}u-1+l_{1}+\de_{i_{1}}(r)+r/T}}{z^{l_{3}+i+1}}Y(u,z)v.
$$

If $g=1$, then $\ast_{g,m,p}^{n}$ is  the same as $\ast_{m,p}^{n}$ defined in
\cite{DJ1}. If $m=p=n$, $\ast_{g,m,p}^{n}$ is just $*_{g,n}$
which has been defined in \cite{DLM4} (see Section 2 of this paper). As in \cite{DLM4}, we will
denote the product by $*_{g,n}$ in this paper.

 Let
$O'_{g,n,m}(V)$ be the linear span of $u\circ_{g,m}^{n}v$ and
 $(L(-1)+L(0)+m-n)u$, where for homogeneous $u\in V^{r}$ and $v\in V$,
$$
u\circ_{g,m}^{n}v={\rm
Res}_{z}\frac{(1+z)^{{\wt}u-1+\de_{i_{1}}(r)+l_{1}+r/T}}{z^{l_{1}+l_{3}+\de_{i_{1}}(r)+\de_{i_{3}}(T-r)+1}}Y(u,z)v.$$
Again if $m=n,$ $u\circ_{g,m}^{n}v=u\circ_{g,n} v$ has been
defined in Section 2 (see also \cite{DLM4}). So if $m=n$,
$O'_{g,n,m}(V)=O_{g,n}(V)$.

\begin{lem}\label{l3.1} If $\overline{i_{1}-i_{3}}\neq r$, then $V^{r}\subseteq
O'_{g,n,m}(V)$.
\end{lem}

\pf Let $u\in V^r$ be homogeneous, then $u\circ_{g,m}^{n}{\bf
1}\in O'_{g,n,m}(V)$. By the definition of $\circ_{g,m}^{n}$, we
have
\begin{eqnarray*}
& & \ \ \ u\circ_{g,m}^{n}{\bf 1}\\
& &
=\sum\limits_{j=0}^{\infty}{{\wt}u-1+\de_{i_{1}}(r)+l_{1}+r/T\choose
j}u_{j-l_{1}-l_{3}-\de_{i_{1}}(r)-\de_{i_{3}}(T-r)-1}{\bf 1}\\
&
&=\sum\limits_{j=0}^{l_{1}+l_{3}+\de_{i_{1}}(r)+\de_{i_{3}}(T-r)}{{\wt}u-1+\de_{i_{1}}(r)+l_{1}+r/T\choose
j}u_{j-l_{1}-l_{3}-\de_{i_{1}}(r)-\de_{i_{3}}(T-r)-1}{\bf 1}.
\end{eqnarray*}
Using relations $u_{-s-1}{\bf 1}=(1/s!)L(-1)^{s}u$ for $s\geq 0$ and
$L(-1)u\equiv (-L(0)-m+n)u$ modulo $O'_{g,n,m}(V)$, we have
\begin{eqnarray*}
& & \ \ \ u\circ_{g,m}^{n}{\bf 1}\\
& &=\sum\limits_{j=0}^{k}{{\wt}u-1+\de_{i_{1}}(r)+l_{1}+r/T\choose
j}(-1)^{k-j}\\
& & \ \ \
\cdot{{\wt}u+2l_{1}+\de_{i_{1}}(r)+\de_{i_{3}}(T-r)-j-1+i_{1}/T-i_{3}/T\choose
k-j}u\\
&
&=\frac{1}{k!}((r/T-(i_{1}-i_{3})/T)^{k}+\sum\limits_{\stackrel{i,j\in\Z_{+}}{0\leq
i+j<k}}a_{i,j}(r/T)^{i}(i_{1}/T-i_{3}/T)^{j})u,
\end{eqnarray*}
where $k=l_{1}+l_{3}+\de_{i_{1}}(r)+\de_{i_{3}}(T-r)$ and
$a_{i,j}\in\Z$. By the fact that $\overline{i_{1}-i_{3}}\neq r$,
we know that $u\circ_{g,m}^{n}{\bf 1}\equiv cu$ modulo
$O'_{g,n,m}(V)$ for a non-zero constant $c$. This shows $u\in
O'_{g,n,m}(V)$. \qed

\begin{coro}\label{co3.2} Let $u\in V^r, v\in V^s$ be homogeneous.
If $\overline{i_{1}-i_{2}}\neq s$, then $u*_{g,m,p}^{n}v\in
O'_{g,n,m}(V)$.
\end{coro}

\pf If $\overline{i_2-i_3}\ne r$ then  $u*_{g,m,p}^{n}v=0$ by
definition. If $\overline{i_2-i_3}=r$ then $u*_{g,m,p}^{n}v\in
V^{r+s}$ and $\overline{i_1-i_3}\ne\overline{r+s}.$ The corollary
follows from Lemma \ref{l3.1}. \qed

The proof of the following lemma is fairly standard (cf.
\cite{DLM3} and \cite{Z}).
\begin{lem}\label{l3.2} For homogeneous  $u,v\in V$, and integers
 $k\geq s\geq 0$,
 $$
{\rm
Res}_{z}\frac{(1+z)^{{\wt}u-1+\de_{i_{1}}(r)+l_{1}+r/T+s}}{z^{l_{1}+l_{3}+\de_{i_{1}}(r)+
\de_{i_{3}}(T-r)+1+k}}Y(u,z)v\in
O'_{g,n,m}(V).$$
\end{lem}

\begin{lem}\label{l3.3}
For homogeneous $u\in V^r$ and $v\in V^s$, if
$\overline{i_{2}-i_{3}}=r$, $\overline{i_{1}-i_{2}}= s$, and
$m+n-p\geq 0$, then
$$u{\ast}_{g,m,p}^{n}v-v\ast_{g,m,m+n-p}^{n}u-{\rm Res}_{z}(1+z)^{{\wt}u-1+p-n}Y(u,z)v\in
O'_{g,n,m}(V).$$
\end{lem}

\pf From the assumption that  $\overline{i_{2}-i_{3}}=r$ and
$\overline{i_{1}-i_{2}}=s$, one can easily deduce that
$-1+\de_{i_{1}}(s)+\de_{i_{2}}(T-s)=0$ and
$-1+\de_{i_{1}}(r)+\de_{i_{3}}(T-r)=\varepsilon$, where
\begin{equation}\label{ea3.0}
\varepsilon=\left\{\begin{array}{lr}1 \quad {\rm if} \
i_{1}+i_{3}-i_{2}\geq T, \\  0 \quad {\rm if} \ 0\leq
i_{1}+i_{3}-i_{2}< T, \\ -1 \quad {\rm if} \ i_{1}+i_{3}-i_{2}< 0.
\end{array}\right.
\end{equation}
From the definition of $O'_{g,n,m}(V)$, we have
$$Y(v,z)u\equiv(1+z)^{-{\wt}u-{\wt}v-m+n}Y(u,\frac{-z}{1+z})v$$
modulo $O'_{g,n,m}(V)$ (cf. \cite{Z} and \cite{DLM2}). Hence
\begin{eqnarray*}
& &\ \ \ \ \
v\ast_{g,m,m+n-p}^{n}u\\
& & =\sum\limits_{i=0}^{l_{1}+l_{3}-l_{2}+\varepsilon}
(-1)^{i}{l_{2} +i\choose i} {\rm
Res}_{z}\frac{(1+z)^{{\wt}v-1+l_{1}+\de_{i_{1}}(s)+s/T}}
{z^{l_{2}+i+1}}Y(v,z)u\\
& &\ \ \
\equiv\sum\limits_{i=0}^{l_{1}+l_{3}-l_{2}+\varepsilon}(-1)^{i}{
l_{2}+i\choose i}{\rm
Res}_{z}\frac{(1+z)^{-{\wt}u-1+l_{3}+\de_{i_{1}}(s)+(s-i_{1}+i_{3})/T}}
{z^{l_{2}+i+1}}Y(u,\frac{-z}{1+z})v\\
& & \ \ \ \ \ \ \ \
 \ ({{\rm mod}}\ O'_{g,n,m}(V))\\
 & &\ \ \
=\sum\limits_{i=0}^{l_{1}+l_{3}-l_{2}+\varepsilon}(-1)^{l_{2}}
{l_{2} +i\choose i}{\rm
Res}_{z}\frac{(1+z)^{{\wt}u-1+l_{2}-l_{3}+i+(i_{2}-i_{3})/T}}{z^{l_{2}+i+1}}Y(u,z)v.
\end{eqnarray*}
Recall the definition of $u{\ast}_{g,m,p}^{n}v:$
\begin{eqnarray*}
& &  u*_{g,m,p}^{n}v=\sum\limits_{i=0}^{l_{2}}
(-1)^{i}{l_{1}+l_{3}-l_{2}+\varepsilon+i\choose
i}\\
& & \ \ \ \ \ \ \ \ \ \ \ \ \ \ \ \ \
\cdot\Res_{z}\frac{(1+z)^{{\wt}u-1+l_{1}+\delta_{i_{1}}(r)+r/T}}{z^{l_{1}+l_{3}-l_{2}+
\varepsilon+i+1}}Y(u,z)v.
\end{eqnarray*}
Since $(i_{3}+r-i_{2})/T=\de_{i_{3}}(T-r)$, we have $$
u{\ast}_{g,m,p}^{n}v-v\ast_{g,m,m+n-p}^{n}u\equiv {\rm
Res}_{z}A_{l_{2},l_{1}+l_{3}-l_{2}+\varepsilon}(z)(1+z)^{{\wt}u-1+p-n}Y(u,z)v$$
modulo $O'_{g,n,m}(V)$ where
\begin{eqnarray*}
& &
A_{l_{2},l_{1}+l_{3}-l_{2}+\varepsilon}(z)=\sum\limits_{i=0}^{l_{2}}(-1)^{i}{
l_{1}+l_{3}-l_{2}+\varepsilon+i\choose
i}\frac{(1+z)^{l_{1}+l_{3}-l_{2}+\varepsilon+1}}{z^{l_{1}+l_{3}-l_{2}+\varepsilon+i+1}}\\
& & \ \ \ \ -
\sum\limits_{i=0}^{l_{1}+l_{3}-l_{2}+\varepsilon}(-1)^{l_{2}}{
l_{2}+i\choose i}\frac{(1+z)^i}{z^{l_{2}+i+1}}.
\end{eqnarray*}
 The lemma now follows from Proposition 5.1 of \cite{DJ1}. \qed

By  Lemma \ref{l3.3} and the fact that ${\bf
1}\bar{*}_{g,m}^{n}u=u$, we have
\begin{coro}\label{co3.1} Let $u\in V^{r}$ be homogeneous,
then  $$ u*_{g,m}^{n}{\bf 1}-u\in O'_{g,n,m}(V).$$
\end{coro}

\begin{lem}\label{l3.4} $V\bar*_{g,m}^{n}O'_{g,n,m}(V)\subseteq
O'_{g,n,m}(V)$, $O'_{g,n,m}(V)*_{g,m}^{n}V\subseteq
O'_{g,n,m}(V)$.
\end{lem}

\pf  The proof is similar to that of Lemma 2.5 in \cite{DJ1}.
Let $u\in V^r$, $v\in V^s$, $w\in V^q$ be homogeneous. By the
definition of $\bar*_{g,m}^{n}$, we can assume that $r=0$. So
\begin{eqnarray*}
& &\ \ \ \ \ u\bar{\ast}_{g,m}^{n}(v\circ_{g,m}^{n}w)\\
& & \equiv\sum\limits_{i=0}^{l_{3}} (-1)^{i}{l_{1}+i\choose
i}\Res_{z_{1}}\frac{(1+z_{1})^{{\wt}u+l_{1}}}{z_{1}^{l_{1}+i+1}}Y(u,z_{1})\\
& & \ \ \ \ \ \ {\rm
Res}_{z_{2}}\frac{(1+z_{2})^{{\wt}v-1+\de_{i_{1}}(s)+l_{1}+s/T}}{z_{2}^{l_{1}+l_{3}+\de_{i_{1}}(s)+\de_{i_{3}}(T-s)+1}}
Y(v,z_{2})w\\
& &\ \ \ \ -\sum\limits_{i=0}^{l_{3}} (-1)^{i}{l_{1}+i\choose
i}{\rm
Res}_{z_{2}}\frac{(1+z_{2})^{{\wt}v-1+\de_{i_{1}}(s)+l_{1}+s/T}}{z_{2}^{l_{1}+l_{3}+\de_{i_{1}}(s)+\de_{i_{3}}(T-s)+1}}
Y(v,z_{2})\\
& & \ \ \ \ \ \ {\rm
Res}_{z_{1}}\frac{(1+z_{1})^{{\wt}u+l_{1}}}{z_{1}^{l_{1}+i+1}}Y(u,z_{1})w\\
& &=\sum\limits_{i=0}^{l_{3}} (-1)^{i}{l_{1}+i\choose
i}\sum\limits_{j\geq0}{{\wt}u+l_{1}\choose
j}\sum\limits_{k=0}^{\infty}
{-l_{1}-i-1\choose k}\\
& &\ \ \ \ \cdot {\rm Res}_{z_{2}}{\rm
Res}_{z_{1}-z_{2}}\frac{(1+z_{2})^{{\wt}u+2l_{1}-j+{\wt}v-1+\de_{i_{1}}(s)+s/T}(z_{1}-z_{2})^{j+k}}{z_{2}^{2l_{1}+i+2+k+l_{3}
+\de_{i_{1}}(s)+\de_{i_{3}}(T-s) }}
Y(Y(u,z_{1}-z_{2})v,z_{2})w\\
& &=\sum\limits_{i=0}^{l_{3}} (-1)^{i}{l_{1}+i\choose
i}\sum\limits_{j\geq0}{{\wt}u+l_{1}\choose
j}\sum\limits_{k=0}^{\infty}
{-l_{1}-i-1\choose k}\\
& &\ \ \ \ \cdot {\rm
Res}_{z_{2}}\frac{(1+z_{2})^{{\wt}u+2l_{1}-j+{\wt}v-1+\de_{i_{1}}(s)+s/T}}{z_{2}^{2l_{1}+i+2+k+l_{3}
+\de_{i_{1}}(s)+\de_{i_{3}}(T-s) }}Y(u_{j+k}v,z_{2})w.
\end{eqnarray*}
Note that the weight of $u_{j+k}v$ is ${\wt}u+{\wt}v-j-k-1.$ By
Lemma \ref{l3.2} we see that
$u\bar{\ast}_{g,m}^{n}(v\circ_{g,m}^{n}w)$ lies in
$O'_{g,n,m}(V).$

By the definition of $*_{g,m}^{n}$ and Corollary \ref{co3.2},
$(v\circ_{g,m}^{n}w){\ast}_{g,m}^{n}u\in O'_{g,n,m}(V)$, if
$\overline{i_{1}-i_{3}}\neq \overline{s+q}$ or $r\neq 0$. So we
can assume that $\overline{i_{1}-i_{3}}= \overline{s+q}$ and
$r=0$. By Lemma \ref{l3.3}, we have
\begin{eqnarray*}
& &\ \ \ \ u\bar{\ast}_{g,m}^{n}(v\circ_{g,m}^{n}w)-(v\circ_{g,m}^{n}w){\ast}_{g,m}^{n}u\\
& &\equiv{\rm Res}_{z}(1+z)^{{\wt}u-1}Y(u,z)(v\circ_{g,m}^{n}w)\\
& &={\rm Res}_{z_{1}}(1+z_{1})^{{\wt}u-1}{\rm
Res}_{z_{2}}\frac{(1+z_{2})^{{\wt}v-1+\de_{i_{1}}(s)+l_{1}+s/T}}{z_{2}^{l_{1}+l_{3}+
\de_{i_{1}}(s)+\de_{i_{3}}(T-s)+1}}Y(u,z_{1})Y(v,z_{2})w\\
& &\equiv\sum\limits_{j\geq 0}{ {\wt}u-1\choose j}{\rm
Res}_{z_{2}}{\rm
Res}_{z_{1}-z_{2}}\frac{(1+z_{2})^{{\wt}u-j+{\wt}v-2+\de_{i_{1}}(s)+l_{1}+s/T}}{z_{2}^{l_{1}+l_{3}+
\de_{i_{1}}(s)+\de_{i_{3}}(T-s)+1}}\\
& &\ \ \ \ \cdot (z_{1}-z_{2})^{j}Y(Y(u,z_{1}-z_{2})v,z_{2})w\\
& &=\sum\limits_{j\geq 0}{ {\wt}u-1\choose j}{\rm
Res}_{z_{2}}\frac{(1+z_{2})^{{\wt}u-j+{\wt}v-2+\de_{i_{1}}(s)+l_{1}+s/T}}{z_{2}^{l_{1}+l_{3}+
\de_{i_{1}}(s)+\de_{i_{3}}(T-s)+1}}Y(u_jv,z_{2})w\in O'_{g,n,m}(V)
\end{eqnarray*}
This proves that $(v\circ_{g,m}^{n}w){\ast}_{g,m}^{n}u\in
O'_{g,n,m}(V).$

Finally we deal with $(L(-1)u+(L(0)+m-n)u)\ast_{g,m}^{n}v$ and
$v\bar{\ast}_{g,m}^{n}(L(-1)u+(L(0)+m-n)u)$.  As before we assume
that $u\in V^r, v\in V^s$ are homogeneous,
$\overline{i_{1}-i_{3}}=r$ and $s=0$. Then
$-1+\de_{i_{1}}(r)+\de_{i_{3}}(T-r)=0$. So
\begin{eqnarray*}
& &\ \ \ \ (L(-1)u+(L(0)+m-n)u){\ast}_{g,m}^{n}v\\
& &=\sum\limits_{i=0}^{l_{1}} (-1)^{i}{l_{3}+i\choose
i}\Res_{z}\frac{(1+z)^{{\wt}u+l_{1}+\de_{i_{1}}(r)+r/T}}{z^{l_{3}+i+1}}Y(L(-1)u,z)v\\
& &\ \ \ \ + ({\wt}u+m-n)\sum\limits_{i=0}^{l_{1}}
(-1)^{i}{l_{3}+i\choose
i}\Res_{z}\frac{(1+z)^{{\wt}u-1+l_{1}+\de_{i_{1}}(r)+r/T}}{z^{l_{3}+i+1}}Y(u,z)v\\
& &=\sum\limits_{i=0}^{l_{1}} (-1)^{i}{l_{3}+i\choose
i}(-l_{3}-\de_{i_{1}}(r)-(r+i_{3}-i_{1})/T)\\
& & \ \ \ \ \ \ \Res_{z}\frac{(1+z)^{{\wt}u-1+l_{1}+\de_{i_{1}}(r)+r/T}}{z^{l_{3}+i+1}}Y(u,z)v\\
& &\ \ \ \ +\sum\limits_{i=0}^{l_{1}}(-1)^{i}{ l_{3}+i\choose i}
(l_{3}+i+1){\rm
Res}_{z}\frac{(1+z)^{{\wt}u+l_{1}+\de_{i_{1}}(r)+r/T}}{z^{l_{3}+i+2}}Y(u,z)v\\
& &=(-1)^{l_{1}}(l_{1}+l_{3}+1){l_{1}+l_{3}\choose l_{1}}{\rm
Res}_{z}\frac{(1+z)^{{\wt}u-1+l_{1}+\de_{i_{1}}(r)+r/T}}{z^{l_{1}+l_{3}+2}}Y(u,z)v.
\end{eqnarray*}
Since $\de_{i_{1}}(r)+\de_{i_{3}}(T-r)+1=2$, it follows that the
last expression is in $O'_{g,n,m}(V)$. Thus
$(L(-1)u+(L(0)+m-n)u){\ast}_{g,m}^{n}v$ belongs to
$O'_{g,n,m}(V).$ We now turn to
$v\bar{\ast}_{g,m}^{n}(L(-1)+L(0)+m-n)u.$ By Lemma \ref{l3.3}, we
have
\begin{eqnarray*}
& &\ \ \ \ v\bar{\ast}_{g,m}^{n}(L(-1)u+(L(0)+m-n)u)-(L(-1)u+(L(0)+m-n)u){\ast}_{g,m}^{n}v\\
& &\equiv{\rm Res}_{z}(1+z)^{{\wt}v-1}Y(v,z)(L(-1)u+(L(0)+m-n)u)\\
& &=\sum\limits_{i\geq0}{ {\wt}v-1\choose
i}v_iL(-1)u+({\wt}u+m-n)\sum\limits_{i\geq 0}{ {\wt}v-1\choose i}v_iu\\
& &=L(-1)\sum\limits_{i\geq 0}{ {\wt}v-1\choose
i}v_iu+\sum\limits_{i\geq 0}{ {\wt}v-1\choose i}iv_{i-1}u\\
& &\ \ \ \ +({\wt}u+m-n)\sum\limits_{i\geq 0}
{ {\wt}v-1\choose i}v_iu\\
& &=L(-1)\sum\limits_{i\geq 0}{ {\wt}v-1\choose i}
v_iu+\sum\limits_{i\geq 0}{ {\wt}v-1\choose i+1}(i+1)v_iu\\
& &\ \ \ \ +({\wt}u+m-n)\sum\limits_{i\geq 0}{ {\wt}v-1\choose i}
v_iu\\
& &=\sum\limits_{i\geq 0} { {\wt}v-1\choose
i}(L(-1)+{\wt}v-i-1+{\wt}u+m-n)v_iu\\
& &= \sum\limits_{i\geq 0} { {\wt}v-1\choose i}(L(-1)v_iu+
L(0)v_iu+ (m-n)v_iu)
\end{eqnarray*}
which is in $O'_{g,n,m}(V).$ So
$v\bar{\ast}_{g,m}^{n}(L(-1)u+(L(0)+m-n)u)\in O'_{g,n,m}(V),$ as
desired. \qed

\begin{lem}\label{l3.5}  We have
$(a\bar{\ast}_{g,m}^{n}b)\ast_{g,m}^{n}c-a\bar{\ast}_{g,m}^{n}(b\ast_{g,m}^{n}c)\in
O'_{g,n,m}(V)$ for homogeneous $a,b,c\in V.$
\end{lem}

\pf Let $a\in V^{r},b\in V^{s},c\in V^l$ be homogeneous. By the
definition of $*_{g,m,p}^{n}$ and Lemma \ref{l3.4},
$(a\bar{\ast}_{g,m}^{n}b)\ast_{g,m}^{n}c$,
$a\bar{\ast}_{g,m}^{n}(b\ast_{g,m}^{n}c)\in O'_{g,n,m}(V)$, if
$r\neq 0$ or $\overline{i_{1}-i_{3}}\neq s$ or $l\neq 0$. So we
can assume that $r=l=0$ and $\overline{i_{1}-i_{3}}=s$. So
\begin{eqnarray*}
& & \ \ \ \ (a\bar{\ast}_{g,m}^{n}b)\ast_{g,m}^{n}c\\
& &=\sum\limits_{k=0}^{l_{1}} (-1)^{k}{l_{3}+k\choose
k}\sum\limits_{i=0}^{l_{3}} (-1)^{i}{l_{1}+i\choose
i}\sum\limits_{j\geq 0}
{{\wt}a+l_{1}\choose j}\\
& &\ \ \ \ \cdot{\rm
Res}_{z}\frac{(1+z)^{{\wt}a+{\wt}b-1+2l_{1}-j+i+\de_{i_{1}}(s)+s/T}}{z^{l_{3}+k+1}}Y(a_{j-l_{1}-i-1}b,z)c
\\
& &=\sum\limits_{k=0}^{l_{1}} (-1)^{k}{l_{3}+k\choose
k}\sum\limits_{i=0}^{l_{3}} (-1)^{i}{l_{1}+i\choose i}{\rm
Res}_{z_{2}}{\rm
Res}_{z_{1}-z_{2}}(z_{1}-z_{2})^{-l_{1}-i-1}\\
& &\ \ \ \
\cdot\frac{(1+z_{1})^{{\wt}a+l_{1}}(1+z_{2})^{{\wt}b-1+l_{1}+i+\de_{i_{1}}(s)+s/T}}{z_{2}^{l_{3}+k+1}}
Y(Y(a,z_{1}-z_{2})b,z_{2})c\\
& &=\sum\limits_{k=0}^{l_{1}} (-1)^{k}{l_{3}+k\choose
k}\sum\limits_{i=0}^{l_{3}} (-1)^{i}{l_{1}+i\choose
i}\sum\limits_{j\geq 0}
{-l_{1}-i-1\choose j}\\
& &\ \ \ \ \cdot{\rm Res}_{z_{1}}{\rm Res}_{z_{2}}
\frac{(1+z_{1})^{{\wt}a+l_{1}}(1+z_{2})^{{\wt}b-1+l_{1}+i+\de_{i_{1}}(s)+s/T}
(-z_{2})^{j}}{z_{1}^{l_{1}+i+1+j}z_{2}^{l_{3}+k+1}}Y(a,z_{1})Y(b,z_{2})c\\
& &\ \ \ \ -\sum\limits_{k=0}^{l_{1}} (-1)^{k}{l_{3}+k\choose
k}\sum\limits_{i=0}^{l_{3}} (-1)^{i}{l_{1}+i\choose
i}\sum\limits_{j\geq 0}
{-l_{1}-i-1\choose j}\\
& &\ \ \ \ \cdot{\rm Res}_{z_{2}}{\rm Res}_{z_{1}}
\frac{(1+z_{1})^{{\wt}a+l_{1}}(1+z_{2})^{{\wt}b-1+l_{1}+i+\de_{i_{1}}(s)+s/T}
z_{1}^{j}}{(-z_{2})^{l_{1}+i+1+j}z_{2}^{l_{3}+k+1}}Y(b,z_{2})Y(a,z_{1})c\\\
& & \equiv
a\bar{\ast}_{g,m}^{n}(b\ast_{g,m}^{n}c)+\sum\limits_{k=0}^{l_{1}}(-1)^{k}
{l_{3}+k\choose k}\sum\limits_{i=0}^{l_{3}}(-1)^{i}
{l_{1}+i\choose i}\\
& &\ \ \ \ \cdot{\rm Res}_{z_{1}}{\rm
Res}_{z_{2}}\left[\sum\limits_{j=0}^{l_{3}-i}(-1)^{j}
{-l_{1}-i-1\choose j}\sum\limits_{q=0}^{i}{ i\choose q}\frac{z_{2}^{j+q}}{z_{1}^{j+i}}-\frac{1}{z_{1}^{i}}\right]\\
& &\ \ \ \
\cdot\frac{(1+z_{1})^{{\wt}a+l_{1}}(1+z_{2})^{{\wt}b-1+l_{1}+\de_{i_{1}}(s)+s/T}}
{z_{1}^{l_{1}+1}z_{2}^{l_{3}+k+1}}Y(a,z_{1})Y(b,z_{2})c.
\end{eqnarray*}
The lemma then follows from Proposition 5.2 in \cite{DJ1}. \qed

 Let $O''_{g,n,m}(V)$ be the linear span of
$u\ast_{g,m,p_{3}}^{n}((a\ast_{g,p_{1},p_{2}}^{p_{3}}b){\ast}_{g,m,p_{1}}^{p_{3}}c-
a{\ast}_{g,m,p_{2}}^{p_{3}}(b{\ast}_{g,m,p_{1}}^{p_{2}}c)),$
 for all $a,b,c,u\in V$ and all $p_{1},p_{2},p_{3}\in\frac{1}{T}{\mathbb
Z}_{+}$.   Set
$$O'''_{g,n,m}(V)=\sum_{p_{1},p_{2}\in\frac{1}{T}{\mathbb
Z}_{+}}(V\ast_{g,p_{1},p_{2}}^{n}O'_{g,p_{2},p_{1}}(V))\ast_{g,m,p_{1}}^{n}V,$$
and
$$
O_{g,n,m}(V)=O'_{g,n,m}(V)+O''_{g,n,m}(V)+O'''_{g,n,m}(V).
$$

\begin{lem}\label{l3.6} For
any $m,n,p\in\frac{1}{T}{\mathbb Z}_{+},$ we have
$$V{\ast}_{g,m,p}^{n}O_{g,p,m}(V)\subseteq O_{g,n,m}(V), \
O_{g,n,p}(V)\ast_{g,m,p}^{n}V\subseteq O_{g,n,m}(V).$$In
particular, $V\bar{\ast}_{g,m}^{n}O_{g,n,m}(V)\subseteq
O_{g,n,m}(V)$, $O_{g,n,m}(V)\ast_{g,m}^{n}V\subseteq
O_{g,n,m}(V)$.
\end{lem}

\pf Note that ${\bf 1}\bar*_{g,m}^{n}u=u$ for any $u\in V.$ We have
\begin{equation}\label{equation1}
O'_{g,n,p}(V)\ast_{g,m,p}^{n}V\subseteq
(V\bar{\ast}_{g,p}^{n}O'_{g,n,p}(V))\ast_{g,m,p}^{n}V\subseteq
O'''_{g,n,m}(V)\end{equation} and
\begin{equation}\label{ea3.1}(a\ast_{g,p_{1},p_{2}}^{n}b){\ast}_{g,m,p_{1}}^{n}c-
a{\ast}_{g,m,p_{2}}^{n}(b{\ast}_{g,m,p_{1}}^{p_{2}}c)\in
O''_{g,n,m}(V)
\end{equation}
for $a,b,c\in V$ and $p,p_1,p_2\in\frac{1}{T}\Z_+.$ Thus
\begin{equation}\label{ea3.2}
V\ast_{g,m,p_{2}}^{n}(O'_{g,p_{2},p_{1}}(V)\ast_{g,m,p_{1}}^{p_{2}}V)\subseteq
O_{g,n,m}(V).
\end{equation} Using Corollary \ref{co3.1} and the definition of $O'''_{g,n,m}(V)$ gives
\begin{equation}\label{equation2}V{\ast}_{g,m,p}^{n}O'_{g,p,m}(V)\subseteq
(V\ast_{g,m,p}^{n}O'_{g,p,m}(V))\ast_{g,m}^{n}V+O'_{g,n,m}(V)\subseteq
O_{g,n,m}(V).
\end{equation}
By (\ref{equation1}) and (\ref{equation2}), it is enough to prove
that
\begin{equation}\label{ea3.3}
V{\ast}_{g,m,p}^{n}((V\ast_{g,p_{1},p_{2}}^{p}O'_{g,p_{2},p_{1}}(V))\ast_{g,m,p_{1}}^{p}V+O''_{g,p,m}(V))\subseteq
O_{g,n,m}(V)
\end{equation}
and
\begin{equation}\label{ea3.4}
((V\ast_{g,p_{1},p_{2}}^{n}O'_{g,p_{2},p_{1}}(V))\ast_{g,p,p_{1}}^{n}V+O''_{g,n,p}(V))\ast_{g,m,p}^{n}V\subseteq
O_{g,n,m}(V)
\end{equation}
for $p_1,p_{2},p\in\frac{1}{T}\Z_+$.

We only prove  (\ref{ea3.3}). The proof of (\ref{ea3.4}) is
similar.  By the definition of $O_{g,n,m}(V)$ and (\ref{ea3.1})
and (\ref{ea3.2}), we
 have
\begin{eqnarray*}
&
&V{\ast}_{g,m,p}^{n}((V\ast_{g,p_{1},p_{2}}^{p}O'_{g,p_{2},p_{1}}(V))\ast_{g,m,p_{1}}^{p}V)\\
& & \subseteq
V{\ast}_{g,m,p}^{n}(V\ast_{g,m,p_{2}}^{p}(O'_{g,p_{2},p_{1}}(V){\ast}_{g,m,p_{1}}^{p_{2}}V))+O_{g,n,m}(V)
 \\
& &\subseteq
(V{\ast}_{g,p_{2},p}^{n}V)\ast_{g,m,p_{2}}^{n}(O'_{g,p_{2},p_{1}}(V){\ast}_{g,m,p_{1}}^{p_{2}}V)+O_{g,n,m}(V)\\
& &\subseteq
O_{g,n,m}(V).
\end{eqnarray*}
It remains to prove that
$V{\ast}_{g,m,p}^{n}O''_{g,p,m}(V)\subseteq O_{g,n,m}(V)$. By
(\ref{ea3.1}), we have
\begin{eqnarray*}
& & \ \ \ \ v{\ast}_{g,m,p}^{n}
(u\ast_{g,m,p_{3}}^{p}((a\ast_{g,p_{1},p_{2}}^{p_{3}}b){\ast}_{g,m,p_{1}}^{p_{3}}c-
a{\ast}_{g,m,p_{2}}^{p_{3}}(b{\ast}_{g,m,p_{1}}^{p_{2}}c)))\\
& &\equiv
(v{\ast}_{g,p_{3},p}^{n}u)\ast_{g,m,p_{3}}^{n}((a\ast_{g,p_{1},p_{2}}^{p_{3}}b){\ast}_{g,m,p_{1}}^{p_{3}}c-
a{\ast}_{g,m,p_{2}}^{p_{3}}(b{\ast}_{g,m,p_{1}}^{p_{2}}c))\\& &
\equiv 0 \ ({\rm mod}\  O_{g,n,m}(V)).
\end{eqnarray*}
\qed

We now define
$$A_{g,n,m}(V)=V/O_{g,n,m}(V).$$
The reason for this definition will become clear from the
$g$-twisted representation theory of $V$ discussed later. If $g=1,$ the
$A_{g,n,m}(V)=A_{n,m}(V)$ has been defined and studied in \cite{DJ1}.

We have  the first  main theorem in this paper.

\begin{theorem}\label{t3.7}  $A_{g,n,m}(V)$ is an $A_{g,n}(V)$-$A_{g,m}(V)$-bimodule such that the
left and right actions of $A_{g,n}(V)$ and $A_{g,m}(V)$ are given
by $\bar*_{g,m}^n$ and $*^n_{g,m}.$
\end{theorem}

\section{ Properties of $A_{g,n,m}(V)$}
\def\theequation{4.\arabic{equation}}
\setcounter{equation}{0}

We will discuss some important properties of $A_{g,n,m}(V)$ in
this section. As in \cite{DJ1}, we will interpret these properties
in terms of twisted representation theory in later sections. In
fact, the twisted representation theory is the origin of the
bimodule $A_{g,n,m}(V)$ and its properties.

First we give an isomorphism between $A_{g,n,m}(V)$ and
$A_{g^{-1},m,n}(V)$ as $A_{g,n}(V)$-$A_{g,m}(V)$-bimodules. 
So we need to define  actions of $A_{g,n}(V)$ and
$A_{g,m}(V)$ on $A_{g^{-1},m,n}(V)$ such that $A_{g^{-1},m,n}(V)$
becomes an $A_{g,n}(V)$-$A_{g,m}(V)$-bimodule. Recall from
\cite{Z} the linear map $\phi:V\to V$ such that
$\phi(v)=e^{L(1)}(-1)^{L(0)}v$ for $v\in V.$ Then from Theorem
\ref{t2.4}, $\phi$ induces an anti-isomorphism from $A_{g,n}(V)$
to $A_{g^{-1},n}(V)$.

\begin{lem}\label{l4.1} The $A_{g^{-1},m,n}(V)$ is an
$A_{g,n}(V)$-$A_{g,m}(V)$-bimodule with
the left action
$\bar\cdot_{g,m}^n$ of $A_{g,n}(V)$ and the right action
$\cdot_{g,m}^n$ of $A_{g,m}(V)$ defined by
 $$u\bar\cdot_{g,m}^nv=v*_{g^{-1},n}^m\phi(u), \ \ v\cdot_{g,m}^n w=\phi(w)\bar *^m_{g^{-1},n}v$$
for $v\in A_{g^{-1},m,n}(V)$, $u\in A_{g,n}(V)$, $w\in
A_{g,m}(V)$.
\end{lem}

\pf The proof is similar to that of Lemma 3.1 of \cite{DJ1}.\qed

\begin{prop}\label{p4.2} The linear map $\phi: A_{g,n,m}(V)\rightarrow
A_{g^{-1},m,n}(V)$ defined by
$$
\phi(u)=e^{L(1)}(-1)^{L(0)}u,
$$
for $u\in A_{g,n,m}(V)$, is an $A_{g,n}(V)-A_{g,m}(V)$-bimodule
isomorphism from $A_{g,n,m}(V)$ to $A_{g^{-1},m,n}(V),$ where the
actions of $A_{g,n}(V)$ and $A_{g,m}(V)$ on $A_{g,n,m}(V)$ are
defined as in Theorem \ref{t3.7}, and the actions on
$A_{g^{-1},m,n}(V)$ are defined as in Lemma \ref{l4.1}.
\end{prop}

\pf Let $m=l_{1}+(i_{1}/T), p=l_{2}+(i_{2}/T),
n=l_{3}+(i_{3}/T)\in\frac{1}{T}\Z$ with $l_{1},l_{2},l_{3}$ three
nonnegative integers and $0\leq i_{1},i_{2},i_{3}\leq T-1$.

We first prove that
\begin{equation}\label{equation3}\phi(O'_{g,n,m}(V))\subset
O'_{g^{-1},m,n}(V). \end{equation} Recall the identities
$$(-1)^{L(0)}Y(u,z)(-1)^{L(0)}=Y((-1)^{L(0)}u,-z)$$
$$e^{L(1)}Y(u,z)e^{-L(1)}=Y(e^{(1-z)L(1)}(1-z)^{-2L(0)}u,\frac{z}{1-z})$$
from \cite{FLM2}.

Take $v={\rm
Res}_{z}\frac{(1+z)^{{\wt}a-1+\de_{i_{1}}(r)+l_{1}+r/T}}{z^{l_{1}+l_{3}+\de_{i_{1}}(r)+\de_{i_{3}}(T-r)+1}}Y(a,z)b\in
O'_{g,n,m}(V)$, where $a\in V^{r}$, $b\in V$. Then
\begin{eqnarray*}
& &\phi(v)={\rm
Res}_{z}e^{L(1)}\frac{(1+z)^{{\wt}a-1+\de_{i_{1}}(r)+l_{1}+r/T}}
{z^{l_{1}+l_{3}+\de_{i_{1}}(r)+\de_{i_{3}}(T-r)+1}}\\
& & \ \ \ \cdot Y(e^{(1+z)L(1)}(1+z)^{-2L(0)}
(-1)^{L(0)}a,\frac{-z}{1+z})e^{L(1)}(-1)^{L(0)}b\\
& &={\rm
Res}_{z}(-1)^{{\wt}a+l_{1}+l_{3}+\de_{i_{1}}(r)+\de_{i_{3}}(T-r)+1}\frac{(1+z)^{{\wt}a+l_{3}+\de_{i_{3}}(T-r)-r/T}}{
z^{l_{1}+l_{3}+\de_{i_{1}}(r)+\de_{i_{3}}(T-r)+1}}\\
& & \ \ \ \ \cdot Y(e^{\frac{1}{{1+z}}L(1)}a,z)e^{L(1)}(-1)^{L(0)}b\\
& &={\rm
Res}_{z}(-1)^{{\wt}a+l_{1}+l_{3}+\de_{i_{1}}(r)+\de_{i_{3}}(T-r)+1}\sum\limits_{j=0}^{\infty}
\frac{1}{j!}\frac{(1+z)^{{\wt}a-j-1+l_{3}+\de_{i_{3}}(T-r)+(T-r)/T}}
{z^{l_{1}+l_{3}+\de_{i_{1}}(r)+\de_{i_{3}}(T-r)+1}}
\\ & & \ \ \ \ \cdot Y(L(1)^{j}a,z)e^{L(1)}(-1)^{L(0)}b
\end{eqnarray*}
which is clearly in $O'_{g^{-1},m,n}(V)$ by the definition of
$O'_{g^{-1},m,n}(V)$. For $u\in V^r$,
\begin{eqnarray*}
& &\ \ \ \ \phi(L(-1)u+(L(0)+m-n)u)\\
& &=e^{L(1)}(-1)^{L(0)}{\rm
Res}_{z}(Y(\omega,z)u+zY(\omega,z)u)+(m-n)e^{L(1)}(-1)^{L(0)}u\\
& &=e^{L(1)}{\rm
Res}_{z}(Y(\omega,-z)+zY(\omega,-z))(-1)^{L(0)}u+(m-n)e^{L(1)}(-1)^{L(0)}u\\
& &={\rm
Res}_{z}(1+z)Y(e^{(1+z)L(1)}(1+z)^{-2L(0)}(-1)^{L(0)}\omega,{\frac{-z}{1+z}})e^{L(1)}(-1)^{L(0)}u\\
& &\ \ \ \ +(m-n)e^{L(1)}(-1)^{L(0)}u \\
& &={\rm Res}_{z}(-(1+z)^{2}+z(1+z))Y(e^{(1+z)^{-1}L(1)}\omega,z)e^{L(1)}(-1)^{L(0)}u\\
& &\ \ \ \ +(m-n)e^{L(1)}(-1)^{L(0)}u\\
& &=-(L(-1)+L(0))e^{L(1)}(-1)^{L(0)}u-(n-m)e^{L(1)}(-1)^{L(0)}u
\end{eqnarray*}
which lies in $O'_{g^{-1},m,n}(V).$
So $\phi(O'_{g,n,m}(V))\subset O'_{g^{-1},m,n}(V)$.

 We next prove that
 \begin{equation}\label{ea4.1}
 \phi(u{\ast}_{g,m,p}^{n}v)\equiv \phi(v){\ast}_{g^{-1},n,p}^{m}\phi(u)
 \end{equation}
modulo $O'_{g^{-1},m,n}(V)$ for $u\in V^{r}, v\in V^{s}.$

 If $\overline{i_{2}-i_{3}}\neq r$, then $\overline{i_{3}-i_{2}}\neq
T-r$. By the definition of $*_{g,m,p}^{n}$ and Corollary
\ref{co3.2}, we have  $u{\ast}_{g,m,p}^{n}v=0$ and
$\phi(v){\ast}_{g^{-1},n,p}^{m}\phi(u)\in O'_{g^{-1},m,n}(V)$, so
(\ref{ea4.1}) holds. Similarly, if $\overline{i_{1}-i_{2}}\neq s$,
 then
$u{\ast}_{g,m,p}^{n}v\in O'_{g,n,m}(V)$ and
$\phi(v){\ast}_{g^{-1},n,p}^{m}\phi(u)=0$. In this case
(\ref{ea4.1}) follows from (\ref{equation3}). So we assume that
$\overline{i_{2}-i_{3}}= r$, $\overline{i_{1}-i_{2}}=s$. Then
$-1+\de_{i_{1}}(r)+\de_{i_{3}}(T-r)=\varepsilon$, where
$\varepsilon$ is defined as in (\ref{ea3.0}). We have the
following computation:
\begin{eqnarray*}
&&\phi(u{\ast}_{g,m,p}^{n}v)=\phi(\sum\limits_{i=0}^{l_{2}}
(-1)^{i}{l_{1}+l_{3}-l_{2}-1+\de_{i_{1}}(r)+\de_{i_{3}}(T-r)+i\choose
i}\\
& & \ \ \ \ \ \ \ \ \ \ \ \ \ \ \ \ \
\cdot\Res_{z}\frac{(1+z)^{{\wt}u-1+l_{1}+\delta_{i_{1}}(r)+r/T}}{z^{l_{1}+l_{3}-l_{2}+
\de_{i_{1}}(r)+\de_{i_{3}}(T-r)+i}}Y(u,z)v)\\
& &=\sum\limits_{i=0}^{l_{2}}
(-1)^{i}{l_{1}+l_{3}-l_{2}+\varepsilon+i\choose
i}\\
& & \ \ \
\cdot\Res_{z}\frac{(1+z)^{{\wt}u-1+l_{1}+\delta_{i_{1}}(r)+r/T}}{z^{l_{1}+l_{3}-l_{2}+\varepsilon
+i+1}}e^{L(1)}Y((-1)^{L(0)}u,-z)(-1)^{L(0)}v\\
& &=\sum\limits_{i=0}^{l_{2}}
(-1)^{i}{l_{1}+l_{3}-l_{2}+\varepsilon+i\choose
i}\Res_{z}\frac{(1+z)^{{\wt}u-1+l_{1}+\delta_{i_{1}}(r)+r/T}}{z^{l_{1}+l_{3}-l_{2}+
\varepsilon+i+1}}\\
& &\ \ \ \ \cdot
Y(e^{(1+z)L(1)}(1+z)^{-2L(0)}(-1)^{L(0)}u,\frac{-z}{1+z})e^{L(1)}(-1)^{L(0)}v
\\
& &=\sum\limits_{i=0}^{l_{2}}
(-1)^{i+l_{1}+l_{3}-l_{2}+\varepsilon+{\wt}u}{l_{1}+l_{3}-l_{2}+\varepsilon+i\choose
i}\\
& & \ \ \ \cdot{\rm
Res}_{z}\frac{(1+z)^{{\wt}u+l_{3}-l_{2}+i+\varepsilon-\de_{i_{1}}(r)-r/T}}{
z^{l_{1}+l_{3}-l_{2}+
\varepsilon+i+1}} Y(e^{(1+z)^{-1}L(1)}u,z)e^{L(1)}(-1)^{L(0)}v\\
& &=\sum\limits_{j=0}^{\infty}\frac{1}{
j!}\sum\limits_{i=0}^{l_{2}}
(-1)^{l_{1}+l_{3}-l_{2}+\varepsilon+{\wt}u}{l_{1}+l_{3}-l_{2}+\varepsilon+i\choose
i}\\
& & \ \ \ \cdot{\rm
Res}_{z}\frac{(1+z)^{{\wt}u+l_{3}-l_{2}+i+\varepsilon-\de_{i_{1}}(r)-r/T-j}}{
z^{l_{1}+l_{3}-l_{2}+
\varepsilon+i+1}}Y(L(1)^{j}u,z)e^{L(1)}(-1)^{L(0)}v\\
& &=\sum\limits_{j=0}^{\infty}\frac{(-1)^{{\wt}u}}{
j!}\sum\limits_{i=0}^{l_{1}+l_{3}-l_{2}+\varepsilon}(-1)^{i}{
l_{2}+i\choose i}\\
& &\ \ \ \ \cdot{\rm
Res}_{z}\frac{(1+z)^{{\wt}u-1-j+l_{3}+\de_{i_{3}}(T-r)+(T-r)/T}}{
z^{l_{2}+i+1}} Y(L(1)^{j}u,z)e^{L(1)}(-1)^{L(0)}v\\
& &\ \ \ \
-\sum\limits_{j=0}^{\infty}\frac{(-1)^{{\wt}u}}{j!}{\rm
Res}_{z}(1+z)^{{\wt}u-j-1+n-p}Y(L(1)^{j}u,z)e^{L(1)}(-1)^{L(0)}v\\
& & \ \ \
 \equiv
\phi(v){\ast}_{g^{-1},n,p}^{m}\phi(u)\ ({\rm mod}\
O'_{g^{-1},m,n}(V)).
\end{eqnarray*}
where we have used Proposition 5.1 of \cite{DJ1} and Lemma
\ref{l3.3} in the last two steps.
 In particular,
$\phi(u\bar{\ast}_{g,m}^{n}v)\equiv
\phi(v)\ast_{g^{-1},n}^{m}\phi(u)$ modulo $O'_{g^{-1},m,n}(V)$ and
$\phi(u{\ast}_{g,m}^{n}v)\equiv
\phi(v)\bar{\ast}_{g^{-1},n}^{m}\phi(u)$ modulo
$O'_{g^{-1},m,n}(V)$ for $u,v\in V.$

As in the proof of Proposition 3.2 in \cite{DJ1}, we can easily
deduce that $\phi(O_{g,n,m}(V))\subseteq O_{g^{-1},m,n}(V)$ by
using (\ref{equation3}) and (\ref{ea4.1}). Thus $\phi:
A_{g,n,m}(V)\to A_{g^{-1},m,n}(V)$ is a well defined bimodule
isomorphism.
 \qed

\begin{prop}\label{p4.3} Let $m,n,l\in\frac{1}{T}\Z_{+}$ such that $m-l,n-l$ are
nonnegative. Then $A_{g,n-l,m-l}(V)$ is an
$A_{g,n}(V)$-$A_{g,m}(V)$-bimodule and the identity map on $V$
induces an epimorphism of $A_{g,n}(V)$-$A_{g,m}(V)$-bimodules from
$A_{g,n,m}(V)$ to $A_{g,n-l,m-l}(V).$
\end{prop}

\pf It is good enough to prove the result for $l=1/T.$ First, from
the definition of $O'_{g,n,m}(V)$ and Lemma \ref{l3.2}, we can
easily see that $O'_{g,n,m}(V)\subseteq O'_{g,n-1/T,m-1/T}(V)$.

Next, we prove
that
\begin{equation}\label{ea4.2}u\ast_{g,p_{1},p_{2}}^{p_{3}}v\equiv
u\ast_{g,p_{1}-1/T,p_{2}-1/T}^{p_{3}-1/T}v \ ({\rm mod} \
O'_{g,p_{3}-1/T,p_{1}-1/T}(V)), \end{equation} for $p_{1}, p_{2},
p_{3}\in\frac{1}{T}{\mathbb Z}_{+}$.

Let $u\in V^r, v\in V^s$ be homogeneous and $p_{i}=s_{i}+j_{i}/T,
i=1,2,3$. We can assume that $\overline {j_{2}-j_{3}}=r$ and
$\overline{j_{1}-j_{2}}=s$.

We first assume  that  $r\neq 0$ and $j_{i}\neq 0, i=1,2,3$.  Then
$j_{3}\neq T-r$. It is easy to see that
$\de_{j_{1}}(r)=\de_{j_{1}-1}(r)+\de_{j_{1},r}$ and
$\de_{j_{3}}(T-r)=\de_{j_{3}-1}(T-r)+\de_{j_{3},T-r}$. So
\begin{eqnarray*}
& &u\ast_{g,p_{1},p_{2}}^{p_{3}}v=\sum\limits_{i=0}^{s_{2}}
(-1)^{i}{s_{1}+s_{3}-s_{2}-1+\de_{j_{1}-1}(r)+\de_{j_{1},r}+\de_{j_{3}-1}(T-r)+i\choose
i}\\
& & \ \ \ \ \ \ \ \ \ \ \ \ \ \ \ \ \
\cdot\Res_{z}\frac{(1+z)^{{\wt}u-1+s_{1}+\delta_{j_{1}-1}(r)+\de_{j_{1},r}+r/T}}{z^{s_{1}+s_{3}-s_{2}+
\de_{j_{1}-1}(r)+\de_{j_{1},r}+\de_{j_{3}-1}(T-r)+i}}Y(u,z)v.
\end{eqnarray*}

If $j_{1}\neq r$, it is clear that
$u\ast_{g,p_{1},p_{2}}^{p_{3}}v=u\ast_{g,p_{1}-1/T,p_{2}-1/T}^{p_{3}-1/T}v$.
If $j_{1}= r$, then
\begin{eqnarray*}
& &u\ast_{g,p_{1},p_{2}}^{p_{3}}v=\sum\limits_{i=0}^{s_{2}}
(-1)^{i}{s_{1}+s_{3}-s_{2}+\de_{j_{3}-1}(T-r)+i\choose
i}\\
& & \ \ \ \ \ \ \ \ \ \ \ \ \ \ \ \ \
\cdot\Res_{z}\frac{(1+z)^{{\wt}u+s_{1}+r/T}}{z^{s_{1}+s_{3}-s_{2}+\de_{j_{3}-1}(T-r)+i+1}}Y(u,z)v\\
 & &= \sum_{i=0}^{s_{2}}{s_{1}+s_{3}-s_{2}+\de_{j_{3}-1}(T-r)+i\choose i}
(-1)^i \Res_z\frac{(1+z)^{{\wt}u+s_{1}+r/T-1}}{z^{s_{1}+s_{3}-s_{2}+\de_{j_{3}-1}(T-r)+i}}Y(u,z)v\\
& & \ \ \ \ \
+\sum_{i=0}^{s_{2}}{s_{1}+s_{3}-s_{2}+\de_{j_{3}-1}(T-r)+i\choose
i}(-1)^i
\Res_z\frac{(1+z)^{{\wt}u+s_{1}+r/T-1}}{z^{s_{1}+s_{3}-s_{2}+\de_{j_{3}-1}(T-r)+i+1}}Y(u,z)v\\
& &\equiv
\sum_{i=0}^{s_{2}}{s_{1}+s_{3}-s_{2}-1+\de_{j_{3}-1}(T-r)+i\choose
i}
(-1)^i \Res_z\frac{(1+z)^{{\wt}u+s_{1}+r/T-1}}{z^{s_{1}+s_{3}-s_{2}+\de_{j_{3}-1}(T-r)+i}}Y(u,z)v\\
& & \ \ \ \ ({\rm mod} \
O'_{g,p_{1}-1/T,p_{3}-1/T}(V))\\
& &=u\ast_{g,p_{1}-1/T,p_{2}-1/T}^{p_{3}-1/T}v.
\end{eqnarray*}

The proof of (\ref{ea4.2}) for other cases is similar. Using
(\ref{ea4.2}), Lemma \ref{l3.6} and the definition of
$O_{g,n,m}(V)$, we conclude that  $O''_{g,n,m}(V)$,
$O'''_{g,n,m}(V)\subset O_{g,n-1/T,m-1/T}(V).$ This together with
 (\ref{ea4.2}) finishes the proof.
\qed

Similar to  Proposition 3.4 in \cite{DJ1}, we have the following
result on tensor products of bimodules.

\begin{prop}\label{p4.4} Define the linear map $\varphi$: $A_{g,n,p}(V)\otimes_{A_{g,p}(V)}A_{g,p,m}(V)\rightarrow
A_{g,n,m}(V)$ by
$$
\varphi(u\otimes v)=u\ast_{g,m,p}^{n}v,$$ for $u\otimes v\in
A_{g,n,p}(V)\otimes_{A_{g,p}(V)}A_{g,p,m}(V)$. Then $\varphi$ is
an $A_{g,n}(V)-A_{g,m}(V)$- bimodule homomorphism from
$A_{g,n,p}(V)\otimes_{A_{g,p}(V)}A_{g,p,m}(V)$ to $A_{g,n,m}(V)$.
\end{prop}

\section{Twisted representation theory}
\def\theequation{5.\arabic{equation}}
\setcounter{equation}{0}

Let $M=\bigoplus_{n\in\frac{1}{T}{\mathbb Z}_{+}}M(n)$ be an admissible
$g$-twisted $V$-module such that $M(0)\ne 0$.
 For homogeneous  $u\in V^{r}$, and $m=l_{1}+i_{1}/T$, $n=l_{2}+i_{2}/T$ such that
 $l_{1},l_{2}\in\Z_{+}$, $0\leq
 i_{1},i_{2}\leq T-1$,
  define the linear map $o_{g,n,m}(u): M(m)\rightarrow M(n)$ by
$$o_{g,n,m}(u)w=u_{{\wt}u+m-n-1}w,$$
 where $w\in M(m)$ and
$u_{{\wt}u+m-n-1}$ is the component operator of
$$Y_M(u,z)=\sum_{n\in
r/T+\Z}u_nz^{-n-1}.$$ Note that if $r\neq\overline{i_{1}-i_{2}}$
then $o_{g,n,m}(u)w=0.$

\begin{lem}\label{l5.1} Let $a\in V^{r}$, $b\in V^{s}$, $m=l_{1}+\frac{i_{1}}{T}$, $p=l_{2}+\frac{i_{2}}{T}$,
$n=l_{3}+\frac{i_{3}}{T}\in\frac{1}{T}\Z_{+}$ with
$l_{1},l_{2},l_{3}$ three nonnegative integers and $0\leq
i_{1},i_{2},i_{3}\leq T-1$. Then
$$
o_{g,n,m}(a\ast_{g,m,p}^{n}b)w=o_{g,n,p}(a)o_{g,p,m}(b)w,$$ where
$w\in M(m)$. In particular,
$$o_{g,n,m}(a\ast_{g,m}^{n}b)w=o_{g,n,m}(a)o_{g,m,m}(b)w, \ \
o_{g,n,m}(a\bar{\ast}_{g,m}^{n}b)w=o_{g,n,n}(a)o_{g,n,m}(b)w.$$
\end{lem}

\pf We first assume that $\overline{i_{1}-i_{2}}\neq s.$ Then
$o_{g,p,m}(b)w=0$. If $\overline{i_{2}-i_{3}}=r$, then
$\overline{i_{1}-i_{3}}\neq \overline{r+s}$ and
$o_{g,n,m}(a\ast_{g,m,p}^{n}b)w=0.$ Otherwise,
$\overline{i_{2}-i_{3}}\neq r$ and $a\ast_{g,m,p}^{n}b=0$. If
$\overline{i_{2}-i_{3}}\neq r$, the proof is similar.

Finally we deal with the case that $\overline{i_{2}-i_{3}}=r$,
$\overline{i_{1}-i_{2}}=s$. Let
$\varepsilon=-1+\de_{i_{1}}(r)+\de_{i_{3}}(T-r)$, then on $M(m)$
we have
\begin{eqnarray*}
& &\ \ \ o_{g,n,m}(a\ast_{g,m,p}^{n}b)\\
& &
 =o_{g,n,m}\left(\sum\limits_{i=0}^{l_{2}}(-1)^{i}{l_{1}+l_{3}-l_{2}+\varepsilon+i \choose i}
{\rm Res}_{z}\frac{(1+z)^{{\wt}a-1+l_{1}+\de_{i_{1}}(r)+r/T}}
{z^{l_{1}+l_{3}-l_{2}+\varepsilon+i+1}}Y(a,z)b\right)\\
 & &=\sum\limits_{i=0}^{l_{2}}(-1)^{i}{l_{1}+l_{3}-l_{2}+\varepsilon+i \choose
i}\sum\limits_{j=0}^{{\wt}a-1+l_{1}+\de_{i_{1}}(r)+r/T}{{\wt}a-1+l_{1}+\de_{i_{1}}(r)+r/T
\choose j}\\& & \ \ \ \
\cdot(a_{j-l_{1}-l_{3}+l_{2}-\varepsilon-i-1}b)_{{\wt}a+{\wt}b-j+2l_{1}-l_{2}+i+\varepsilon+(i_{1}-i_{3})/T-1}
\\
&
&=\sum\limits_{i=0}^{l_{2}}(-1)^{i}{l_{1}+l_{3}-l_{2}+\varepsilon+i
\choose
i}\sum\limits_{j=0}^{{\wt}a-1+l_{1}+\de_{i_{1}}(r)+r/T}{{\wt}a-1+l_{1}+\de_{i_{1}}(r)+r/T
\choose j}\\
& & \ \  \ \cdot{\rm Res}_{z_{2}}{\rm
Res}_{z_{0}}z_{0}^{j-l_{1}-l_{3}+l_{2}-\varepsilon-i-1}z_{2}^{{\wt}a+{\wt}b-j+2l_{1}-l_{2}+i+
\varepsilon+(i_{1}-i_{3})/T-1}
Y_{M}(Y(a,z_{0})b,z_{2})\\
&
&=\sum\limits_{i=0}^{l_{2}}(-1)^{i}{l_{1}+l_{3}-l_{2}+\varepsilon+i
\choose i}{\rm Res}_{z_{2}}{\rm
Res}_{z_{0}}z_{0}^{-l_{1}-l_{3}+l_{2}-\varepsilon-i-1}
\\
& &\ \ \  \ \
\cdot(z_{2}+z_{0})^{{\wt}a-1+l_{1}+\de_{i_{1}}(r)+r/T}
z_{2}^{{\wt}b+l_{1}-l_{2}+i-\de_{i_{1}}(r)-r/T+(i_{1}-i_{3})/T+\varepsilon}
Y_{M}(Y(a,z_{0})b,z_{2})\\
&
&=\sum\limits_{i=0}^{l_{2}}(-1)^{i}{l_{1}+l_{3}-l_{2}+\varepsilon+i
\choose i}\sum\limits_{j=0}^{\infty}{
-l_{1}-l_{3}+l_{2}-\varepsilon-i-1\choose j}{\rm Res}_{z_{1}}{\rm
Res}_{z_{2}}(-z_{2})^{j}\\
& & \ \ \  \ \ \cdot
z_{1}^{{\wt}a-1+l_{1}+\de_{i_{1}}(r)+r/T}z_{1}^{-l_{1}-l_{3}+l_{2}-\varepsilon-i-1-j}\\
& & \ \ \ \  \ \ \cdot
z_{2}^{{\wt}b+l_{1}-l_{2}+i-\de_{i_{1}}(r)-r/T+(i_{1}-i_{3})/T+\varepsilon}Y_{M}(a,z_{1})Y_{M}(b,z_{2})
\\
& &\ \ \
-\sum\limits_{i=0}^{l_{2}}(-1)^{i}{l_{1}+l_{3}-l_{2}+\varepsilon+i
\choose i}\sum\limits_{j=0}^{\infty}{
-l_{1}-l_{3}+l_{2}-\varepsilon-i-1\choose j}\\
& & \ \ \  \ \ {\rm Res}_{z_{2}}{\rm Res}_{z_{1}}z_{1}^{j}
z_{1}^{{\wt}a-1+l_{1}+\de_{i_{1}}(r)+r/T}(-z_{2})^{-l_{1}-l_{3}+l_{2}-\varepsilon-i-1-j}
\\
& & \ \ \  \ \ \cdot
z_{2}^{{\wt}b+l_{1}-l_{2}+i-\de_{i_{1}}(r)-r/T+(i_{1}-i_{3})/T+\varepsilon}
Y_{M}(b,z_{2})Y_{M}(a,z_{1})
\\
&
&=\sum\limits_{i=0}^{l_{2}}(-1)^{i}{l_{1}+l_{3}-l_{2}+\varepsilon+i
\choose i}\sum\limits_{j=0}^{\infty}{
-l_{1}-l_{3}+l_{2}-\varepsilon-i-1\choose j}(-1)^{j}\\
& & \ \ \ \ \ \  \cdot
a_{{\wt}a-2+l_{2}-l_{3}-i-j+\de_{i_{1}}(r)+r/T-\varepsilon}b_{{\wt}b+i+j+l_{1}-l_{2}
-\de_{i_{1}}(r)-r/T+(i_{1}-i_{3})/T+\varepsilon}\\
&  &\ \ \
-\sum\limits_{i=0}^{l_{2}}(-1)^{i}{l_{1}+l_{3}-l_{2}+\varepsilon+i
\choose i}\sum\limits_{j=0}^{\infty}{
-l_{1}-l_{3}+l_{2}-q-i-1\choose j}\\
&  &\ \ \ \ \ \ \cdot(-1)^{-l_{1}-l_{3}+l_{2}-q-i-1-j}
 b_{{\wt}b-l_{3}-j-1-\de_{i_{1}}(r)-r/T+(i_{1}-i_{3})/T}a_{{\wt}a-1+l_{1}+\de_{i_{1}}(r)+r/T+j}\\
 &
 &=a_{{\wt}a-1+l_{2}-l_{3}+i_{2}/T-i_{3}/T}b_{{\wt}b-1+l_{1}-l_{2}+i_{1}/T-i_{2}/T}.
\end{eqnarray*}
The proof is complete. \qed

\begin{lem}\label{l5.2} Let   $m=l_{1}+i_{1}/T,n=l_{3}+i_{3}/T$ such that
 $l_{1},l_{3}\in\Z_{+}$, $0\leq
 i_{1},i_{3}\leq T-1$.  Then
 $$o_{g,n,m}(a)=0$$
on $M(m)$,  for all $a\in O_{g,n,m}(V)$.
\end{lem}

 \pf Let $u\in V^{r}$, $v\in V^{s}$. We first prove that
 $o_{g,n,m}(u\circ_{g,m}^{n}v)=0$. It is clear if $\overline{i_{1}-i_{3}}\neq \overline{r+s}.$ So we
 assume that $\overline{i_{1}-i_{3}}=\overline{r+s}$. Then
 \begin{eqnarray*}
 & & \ \ \ \ \ o_{g,n,m}(u\circ_{g,m}^{n}v)\\
 & & =o_{g,n,m}\left({\rm
Res}_{z}\frac{(1+z)^{{\wt}u-1+\de_{i_{1}}(r)+l_{1}+r/T}}{z^{l_{1}+l_{3}+\de_{i_{1}}(r)+\de_{i_{3}}(T-r)+1}}
Y(u,z)v\right)\\
& &
=o_{g,n,m}\left(\sum\limits_{j=0}^{\infty}{{\wt}u-1+\de_{i_{1}}(r)+l_{1}+r/T
\choose
j}u_{j-l_{1}-l_{3}-\de_{i_{1}}(r)-\de_{i_{3}}(T-r)-1}v\right)\\
& & =\sum\limits_{j=0}^{\infty}{{\wt}u-1+\de_{i_{1}}(r)+l_{1}+r/T
\choose
j}\\
& & \ \ \ \ \ \cdot
(u_{j-l_{1}-l_{3}-\de_{i_{1}}(r)-\de_{i_{3}}(T-r)-1}v)_{{\wt}u+{\wt}v-1-j+2l_{1}+\de_{i_{1}}(r)+\de_{i_{3}}(T-r)
+(i_{1}-i_{3})/T}\\
& & ={\rm Res}_{z_{2}}{\rm
Res}_{z_{0}}\frac{(z_{2}+z_{0})^{{\wt}u-1+\de_{i_{1}}(r)+l_{1}+r/T}z_{2}^{{\wt}v+l_{1}+\de_{i_{3}}(T-r)-r/T
 +(i_{1}-i_{3})/T}}{z_{0}^{l_{1}+l_{3}+\de_{i_{1}}(r)+\de_{i_{3}}(T-r)+1}}\\
 & & \ \ \ \ \ \cdot Y_{M}(Y(u,z_{0}),z_{2})v\\
 &
 &=\sum\limits_{j=0}^{\infty}{-l_{1}-l_{3}-\de_{i_{1}}(r)-\de_{i_{3}}(T-r)-1\choose
 j}(-1)^{j}{\rm Res}_{z_{1}}{\rm
Res}_{z_{2}}\frac{z_{1}^{{\wt}u-1+\de_{i_{1}}(r)+l_{1}+r/T}}
 {z_{1}^{l_{1}+l_{3}+\de_{i_{1}}(r)+\de_{i_{3}}(T-r)+1+j}}\\
 & & \ \ \ \ \ \cdot z_{2}^{{\wt}v+l_{1}+\de_{i_{3}}(T-r)-r/T
 +(i_{1}-i_{3})/T+j}Y_{M}(u,z_{1})Y_{M}(v,z_{2})\\
 & &-\sum\limits_{j=0}^{\infty}{-l_{1}-l_{3}-\de_{i_{1}}(r)-\de_{i_{3}}(T-r)-1\choose
 j}{\rm Res}_{z_{2}}{\rm
Res}_{z_{1}}\frac{z_{2}^{{\wt}v+l_{1}+\de_{i_{3}}(T-r)-r/T
 +(i_{1}-i_{3})/T}}{z_{2}^{l_{1}+l_{3}+\de_{i_{1}}(r)+\de_{i_{3}}(T-r)+1+j}}\\
 & & \ \ \ \ \
 \cdot(-1)^{l_{1}+l_{3}+\de_{i_{1}}(r)+\de_{i_{3}}(T-r)+1+j}z_{1}^{{\wt}u-1+\de_{i_{1}}(r)+l_{1}+r/T+j}
 Y_{M}(v,z_{2})Y_{M}(u,z_{1})\\
 & & =\sum\limits_{j=0}^{\infty}{-l_{1}-l_{3}-\de_{i_{1}}(r)-\de_{i_{3}}(T-r)-1\choose
 j}(-1)^{j}\\
 & & \ \ \  \ \ \cdot u_{{\wt}u-2+r/T-l_{3}-\de_{i_{3}}(T-r)-j}v_{{\wt}v+l_{1}+\de_{i_{3}}(T-r)-r/T
 +(i_{1}-i_{3})/T+j}\\
 & &
 \ \ \ -\sum\limits_{j=0}^{\infty}{-l_{1}-l_{3}-\de_{i_{1}}(r)-\de_{i_{3}}(T-r)-1\choose
 j}(-1)^{l_{1}+l_{3}+\de_{i_{1}}(r)+\de_{i_{3}}(T-r)+1+j}
 \\
 &  &\ \ \  \ \ \cdot v_{{\wt}v-r/T
 +(i_{1}-i_{3})/T-l_{3}-\de_{i_{1}}(r)-1-j}u_{{\wt}u-1+\de_{i_{1}}(r)+l_{1}+r/T+j}\\
 & & =0
\end{eqnarray*}
By Lemma \ref{l5.1} and the definition of $O_{g,n,m}(V)$, we know
that $o_{g,n,m}(a)=0$, for all $a\in
O''_{g,n,m}(V)+O'''_{g,n,m}(V)$. \qed

Let $M=\bigoplus_{m\in\frac{1}{T}\Z_{+}}M(m)$ be an admissible
$g$-twisted $V$-module with $M(0)\ne 0.$ Then ${\rm
Hom}_{\C}(M(m),M(n))$ is an $A_{g,n}(V)$-$A_{g,m}(V)$-bimodule
such that $(a\cdot f\cdot b)(w)=af(bw)$ for $a\in A_{g,n}(V), b\in
A_{g,m}(V),$ $f\in  {\rm Hom}_{\C}(M(m),M(n))$ and $w\in M(m).$
Set
$$o_{g,n,m}(V)=\{o_{g,n,m}(v)|v\in V\}.$$

By Lemmas \ref{l5.1} and \ref{l5.2}, we immediately have

\begin{prop}\label{p5.3}   $o_{g,n,m}(V)$ is an $A_{g,n}(V)$-$A_{g,m}(V)$-subbimodule of
${\rm Hom}_{\C}(M(m),M(n))$ and $v\mapsto o_{g,n,m}(v)$ for $v\in
V$ induces an $A_{g,n}(V)$-$A_{g,m}(V)$-bimodule epimorphism from
$A_{g,n,m}(V)$ to $o_{g,n,m}(V).$
\end{prop}

\begin{prop}\label{p5.4} For any $n\in\frac{1}{T}\Z_{+},$  $A_{g,n}(V)$ and $A_{g,n,n}(V)$ are the same.
\end{prop}

\pf The proof is similar to that of Proposition 4.6 of \cite{DJ1}.  \qed

Next we reconstruct the Verma type  admissible $g$-twisted
$V$-module $\bar{M}(U)$ generated by  an $A_{g,m}(V)$-module $U$
by using the bimodules $A_{g,n,m}(V)$. Note that we do not assume
that $U$ cannot factor through $A_{g,m-1/T}(V)$ at this point.

Set
$$
M(U)=\bigoplus_{n\in\frac{1}{T}{\mathbb
Z}_{+}}A_{g,n,m}(V)\otimes_{A_{g,m}(V)} U.$$
 Then $M(U)$ is  $\frac{1}{T}{\mathbb Z}_{+}$-graded such that
$M(U)(n)=A_{g,n,m}(V)\otimes_{A_{g,m}(V)}U.$ For $u\in V^r,$
$p,n\in\frac{1}{T}\Z$, define an operator $u_{p}$ from $M(U)(n)$
to $M(U)(n+{\wt}u-p-1)$ by
$$
u_{p}(v\otimes
w)=\left\{\begin{array}{rl}(u\ast_{g,m,n}^{{\wt}u-p-1+n}v)\otimes
w, \quad {\rm if} \ {\wt}u-1-p+n\geq 0, \\  0, \quad {\rm if} \
{\wt}u-1-p+n< 0,
\end{array}\right.
$$
 for $v\in A_{g,n,m}(V)$ and $w\in U$. It is clear that if
 $p\notin\Z+r/T$, then $u_{p}=0$.

\begin{lem} The action $u_{p}$ is well defined.
\end{lem}

\pf Let $v\in O_{g,n,m}(V)$ and $w\in U.$ By Lemma \ref{l3.6},
$u\ast_{g,m,n}^{{\wt}u-p-1+n}v\in
V\ast_{g,m,n}^{{\wt}u-p-1+n}O_{g,n,m}(V)\subseteq
O_{g,{\wt}u-p-1+n,m}(V)$, so we have $u_{p}(v\otimes w)=0$. Now
let $a\in A_{g,m}(V)$, $v\in A_{g,n,m}(V)$, $w\in U$. Then
\begin{eqnarray*}
& &u_{p}((v\ast_{g,m}^{n}a)\otimes
w)=(u\ast_{g,m,n}^{{\wt}u-p-1+n}(v\ast_{g,m}^{n}a))\otimes w\\
& &\ \ \ \
=((u\ast_{g,m,n}^{{\wt}u-p-1+n}v)\ast_{g,m}^{{\wt}u-p-1+n}a)\otimes
w\\
& &\ \ \ \ =(u\ast_{g,m,n}^{{\wt}u-p-1+n}v)\otimes a\cdot
w=u_{p}(v\otimes a\cdot w).
\end{eqnarray*}
Thus $u_{p}$ is well defined. \qed

For short we set $M=M(U).$ Also let
$$Y_{M}(u,z)=\sum\limits_{p\in{\mathbb Z}+r/T}u_{p}z^{-p-1},$$
for $u\in V^{r}$. It is our desire to prove that $(M(U), Y_M)$ is
an admissible $g$-twisted
 $V$-module isomorphic to the $\bar{M}(U)$ given in Theorem
\ref{t2.4}.

\begin{lem}\label{l5.5} For homogeneous $u\in V^r$ , $v\otimes w\in A_{g,n,m}(V)\otimes_{A_{g,m}(V)} U$ and
$p\in{\mathbb Z}+r/T$, we have

(1) $u_{p}(v\otimes w)=0$, for $p$ sufficiently large;

(2) $Y_{M}(\1,z)={\rm id}_{M}.$
\end{lem}

\pf From the definition, (1) is obvious. We now prove (2). By the
definition of $u_{p}$, we have
\begin{eqnarray*}
& &\1_{p}(v\otimes w)=(1\ast_{g,m,n}^{-p-1+n}v)\otimes w\\
& &=\sum\limits_{i=0}^{l_{2}} (-1)^{i}{l_{1}-p-1+i\choose
i}\Res_{z}\frac{(1+z)^{l_{1}}}{z^{l_{1}-p+i}}Y(1,z)v\otimes w,
\end{eqnarray*}
where $p\in\Z$, $m=l_{1}+i_{1}/T, n=l_{2}+i_{2}/T$ such that
$l_{1},l_{2}\in\Z_{+}$ and $0\leq i_{1},i_{2}<T$. Thus
$\1_{p}(v\otimes w)=0$ if $p<-1$ and $\1_{-1}(v\otimes w)=v\otimes
w.$ By the definition of ${\bf 1}_{p}$, $\1_{p}(v\otimes w)=0$ if
$p\geq l_{2}$. If $-1<p<l_{2}$, then
\begin{eqnarray*}
 &
&\1_{p}(v\otimes
w)=\sum\limits_{i=0}^{p+1}(-1)^{i}{l_{1}-p-1+i\choose
i}{l_{1}\choose
p+1-i}v\otimes w\\
& &=\sum\limits_{i=0}^{p+1}(-1)^{i+p+1}{l_{1}\choose
p+1}{p+1\choose
i}v\otimes w\\
& &=0.
\end{eqnarray*}
As a result we have $Y_{M}(\1,z)={\rm id}_{M}.$ \qed

The main axiom in the definition of admissible $g$-twisted
$V$-modules is the twisted Jacobi identity. As we have already
mentioned in Section 2, the twisted Jacobi identity is equivalent
to the commutator formula (\ref{cea}) and the associativity
(\ref{ea}).  We have  the commutator formula:
\begin{lem}\label{l5.6} For $a\in V^r$, $b\in V^s$, we have
$$
[Y_{M}(a,z_{1}),Y_{M}(b,z_{2})]={\rm
Res}_{z_{0}}z_{2}^{-1}\left(\frac{z_1-z_0}{z_2}\right)^{-r/T}\delta(\frac{z_{1}-z_{0}}{z_{2}})Y_{M}(Y(a,z_{0})b,z_{2}),$$
or equivalently, for $p\in\Z+r/T$, $q\in\Z+s/T$,
$$[a_p,b_q]=\sum_{i=0}^{\infty}{p\choose i}(a_{i}b)_{p+q-i}.$$
\end{lem}

\pf We need to prove that
$$ (a_{p}b_{q}-b_{q}a_{p})(v\otimes w)=\sum_{i=0}^{\infty}{p\choose i}(a_{i}b)_{p+q-i}(v\otimes w)$$
for $p\in\Z+r/T$, $q\in\Z+s/T$ and $v\otimes w\in
A_{g,n,m}(V)\otimes_{A_{g,m}(V)} U.$ This is clear from the
definition of the actions if ${\wt} a+{\wt} b-p-q-2+n<0.$ We now
assume that ${\wt} a+{\wt} b-p-q-2+n\geq 0.$

If ${\wt}a-p-1+n\geq 0$, ${\wt}b-q-1+n\geq 0$ then by Lemma
\ref{l3.3} we have
\begin{eqnarray*}
& &\ \ \ \  a_{p}b_{q}(v\otimes w)-b_{q}a_{p}(v\otimes w)\\
& &=a_{p}(b\ast_{g,m,n}^{{\wt}b-q-1+n}v)\otimes
w-b_{q}(a\ast_{g,m,n}^{{\wt}a-p-1+n}v)\otimes w\\
&
&=\left(a\ast_{g,m,{\wt}b-q-1+n}^{{\wt}a+{\wt}b-p-q-2+n}(b\ast_{g,m,n}^{{\wt}b-q-1+n}v)\right)\otimes
w\\
& & \ \ \ \ \   -\left(b\ast_{g,m,{\wt}a-p-1+n}^{{\wt}a+{\wt}b-p-q-2+n}(a\ast_{g,m,n}^{{\wt}a-p-1+n}v)\right)\otimes w\\
&
&=\left((a\ast_{g,n,{\wt}b-q-1+n}^{{\wt}a+{\wt}b-p-q-2+n}b)\ast_{g,m,n}^{{\wt}a+{\wt}b-p-q-2+n}v\right)\otimes
w\\
& & \ \ \ \ \
-\left((b\ast_{g,n,{\wt}a-p-1+n}^{{\wt}a+{\wt}b-p-q-2+n}a)\ast_{g,m,n}^{{\wt}a+w
b-p-q-2+n}v\right)\otimes w\\
& &=\left(({\rm
Res}_{z}(1+z)^{p}Y_{M}(a,z)b)\ast_{g,m,n}^{{\wt}a+{\wt}b-p-q-2+n}v\right)\otimes
w\\
& &=(\sum\limits_{i=0}^{\infty}{p\choose
i}\left(a_ib)\ast_{g,m,n}^{{\wt}a+{\wt}b-p-q-2+n}v\right)\otimes
w\\
& &=\sum\limits_{i=0}^{\infty}{p\choose i}(a_ib)_{p+q-i}(v\otimes
w).
\end{eqnarray*}

If ${\wt}a-p-1+n<0$, ${\wt}b-q-1+n\geq 0$ then
$b_{q}a_{p}(v\otimes w)=0$ and
\begin{eqnarray*}
& &a_{p}b_{q}(v\otimes w)-b_{q}a_{p}(v\otimes
w)=\left((a\ast_{g,n,{\wt}b-q-1+n}^{{\wt}a+{\wt}b-p-q-2+n}b)
\ast_{g,m,n}^{{\wt}p+{\wt}b-p-q-2+n}v\right)\otimes
w\\
& &\ \ \ \ =\left(({\rm
Res}_{z}(1+z)^{p}Y_{M}(a,z)b)\ast_{g,m,n}^{{\wt}a+{\wt}b-p-q-2+n}v\right)\otimes
w\\
& &\ \ \ \ =\sum\limits_{i=0}^{\infty}{p\choose
i}(a_ib)_{p+q-i}(v\otimes w)
\end{eqnarray*}
where we have used Lemma \ref{l5.7} below. If ${\wt}a-p-1+n\geq
0$, ${\wt}b-q-1+n<0$ the proof is similar. \qed

\begin{lem}\label{l5.7}
Let $u\in V^r, v\in V^s$ and $m=l_{1}+i_{1}/T, n=l_{3}+i_{3}/T$,
$p=l_{2}+i_{2}/T$ such that $l_{1},l_{3}\in\Z_{+}, l_{2}\in\Z,
0\leq i_{1},i_{2},i_{3}<T$,  and $\overline{i_{2}-i_{3}}=r$,
$\overline{i_{1}-i_{2}}=s$.
 If $p\geq 0$, $m+n-p<0$, then
$$u{\ast}_{g,m,p}^{n}v-{\rm Res}_{z}(1+z)^{{\wt}u-1+p-n}Y(u,z)v\in
O'_{g,n,m}(V)$$ and if $p<0$, $m+n-p\geq 0$, then
$$-v\ast_{g,m,m+n-p}^{n}u-{\rm Res}_{z}(1+z)^{{\wt}u-1+p-n}Y(u,z)v\in
O'_{g,n,m}(V).$$
\end{lem}

\pf We first assume that  $p\geq 0$, $m+n-p<0.$ By (\ref{ea3.0}),
$-l_{1}-l_{3}+l_{2}-\varepsilon-1\in\Z_{+}.$  using the definition
gives
\begin{eqnarray*}
& &  u*_{g,m,p}^{n}v=\sum\limits_{i=0}^{l_{2}}
(-1)^{i}{l_{1}+l_{3}-l_{2}+\varepsilon+i\choose
i}\\
& & \ \ \ \ \ \ \ \ \ \ \ \ \ \ \ \ \
\cdot\Res_{z}\frac{(1+z)^{{\wt}u-1+l_{1}+\delta_{i_{1}}(r)+r/T}}{z^{l_{1}+l_{3}-l_{2}+
\varepsilon+i+1}}Y(u,z)v.
\end{eqnarray*}
 Since
$\Res_zY(u,z)v\frac{(1+z)^{{\wt}u-1+l_{1}+\delta_{i_{1}}(r)+r/T}}{z^{l_{1}+l_{3}-l_{2}+
\varepsilon+i+1}}\in O'_{g,n,m}(V)$ if $i>l_{2}$ by Lemma
\ref{l3.2} we see that
\begin{eqnarray*}
& &u{\ast}_{g,m,p}^{n}v\equiv
\sum_{i=0}^{\infty}{-l_{1}-l_{3}+l_{2}-\varepsilon-1\choose i}
\Res_zY(u,z)v\frac{(1+z)^{{\wt}u-1+l_{1}+\delta_{i_{1}}(r)+r/T}}{z^{l_{1}+l_{3}-l_{2}+
\varepsilon+i+1}}\\
& &\ \ \ \  = \Res_z(1+z)^{{\wt}u-1+p-n}Y(u,z)v,
\end{eqnarray*}
where in the last step we have used the fact that
$\de_{i_{3}}(T-r)=(r+i_{3}-i_{2})/T$. So in this case we are
done.

If  $p<0$, $m+n-p\geq 0$ then the result in the first case gives
$$v\ast_{g,m,m+n-p}^{n}u\equiv \Res_z(1+z)^{{\wt}v-1+m-p}Y(v,z)u$$
modulo $O'_{g,n,m}(V).$ Using the identity
$$Y(v,z)u\equiv(1+z)^{-{\wt}u-{\wt}v-m+n}Y(u,\frac{-z}{1+z})v$$
modulo $O'_{g,n,m}(V)$ we see that
\begin{eqnarray*}
& & \Res_z(1+z)^{{\wt}v-1+m-p}Y(v,z)u\equiv
\Res_z(1+z)^{-{\wt}u-1+n-p}Y(u,\frac{-z}{1+z})v\\
& &\ \ \ \ =-\Res_z(1+z)^{{\wt}u-1+p-n}Y(u,z)v.
\end{eqnarray*}
The proof is complete. \qed

\begin{lem}\label{l5.8} Let $n=l_{3}+i_{3}/T\in\frac{1}{T}\Z_{+}$  with $l_{3}\in\Z_{+}$ and $0\leq i_{3}<T$.
Then for $a\in V^r$ and $i\in{\mathbb Z}_{+},$ we have
\begin{eqnarray*}
& & \ \ \ {\rm
Res}_{z_{0}}z_{0}^{i}(z_{0}+z_{2})^{{\wt}a-1+l_{3}+\de_{i_{3}}(r)+r/T}Y_{M}(a,z_{0}+z_{2})Y_{M}(b,z_{2})
\\
& &={\rm
Res}_{z_{0}}z_{0}^{i}(z_{2}+z_{0})^{{\wt}a-1+l_{3}+\de_{i_{3}}(r)+r/T}Y_{M}(Y(a,z_{0})b,z_{2})
\end{eqnarray*}
on $M(U)(n)$.
\end{lem}

\pf Note that  $a_{{\wt}a-1+l_{3}+\de_{i_{3}}(r)+r/T+j}=0$ on
$M(U)(n)$ for any nonnegative integer $j.$ Then
$$
{\rm
Res}_{z_{1}}(z_{1}-z_{2})^{i}z_{1}^{{\wt}a-1+l_{3}+\de_{i_{3}}(r)+r/T}Y_{M}(b,z_{2})Y_{M}(a,z_{1})=0$$
on $M(U)(n)$ and
\begin{eqnarray*}
& &{\rm
Res}_{z_{0}}z_{0}^{i}(z_{0}+z_{2})^{{\wt}a-1+l_{3}+\de_{i_{3}}(r)+r/T}Y_{M}(a,z_{0}+z_{2})Y_{M}(b,z_{2})\\
& &={\rm
Res}_{z_{1}}(z_{1}-z_{2})^{i}z_{1}^{{\wt}a-1+l_{3}+\de_{i_{3}}(r)+r/T}
(Y_{M}(a,z_{1})Y_{M}(b,z_{2})-Y_{M}(b,z_{2})Y_{M}(a,z_{1}))\\
& &={\rm
Res}_{z_{1}}(z_{1}-z_{2})^{i}z_{1}^{{\wt}a-1+l_{3}+\de_{i_{3}}(r)+r/T}[Y_{M}(a,z_{1}),Y_{M}(b,z_{2})]\\
& &={\rm Res}_{z_{0}}{\rm
Res}_{z_{1}}(z_{1}-z_{2})^{i}z_{1}^{{\wt}a-1+l_{3}+\de_{i_{3}}(r)+r/T}
\\& & \ \ \  \ \ \cdot\left(\frac{z_1-z_0}{z_2}\right)^{-r/T}z_{2}^{-1}
\delta(\frac{z_{1}-z_{0}}{z_{2}})Y_{M}(Y(a,z_{0})b,z_{2})\\
& &={\rm Res}_{z_{0}}{\rm
Res}_{z_{1}}z_{0}^{i}(z_{2}+z_{0})^{{\wt}a-1+l_{3}+\de_{i_{3}}(r)+r/T}\\
& & \ \ \ \ \
\cdot\left(\frac{z_1-z_0}{z_2}\right)^{-r/T}z_{1}^{-1}\delta(\frac{z_{2}+z_{0}}{z_{1}})Y_{M}(Y(a,z_{0})b,z_{2})\\
& &={\rm
Res}_{z_{0}}z_{0}^{i}(z_{2}+z_{0})^{{\wt}a-1+l_{3}+\de_{i_{3}}(r)+r/T}Y_{M}(Y(a,z_{0})b,z_{2}),
\end{eqnarray*}
where we have used Lemma \ref{l5.6}.\qed

\begin{lem}\label{l5.9} Let $n=l_{3}+i_{3}/T\in\frac{1}{T}\Z_{+}$  with $l_{3}\in\Z_{+}$ and $0\leq i_{3}<T$.
Then for $a\in V^r$ and $l\in{\mathbb Z}_{+}\setminus\{0\},$ we
have
\begin{eqnarray*}
& &{\rm
Res}_{z_{0}}z_{0}^{-l}(z_{2}+z_{0})^{{\wt}a+q}z_{2}^{{\wt}b-q}
Y_{M}(Y(a,z_{0})b,z_{2})\\
& &={\rm
Res}_{z_{0}}z_{0}^{-l}(z_{0}+z_{2})^{{\wt}a+q}z_{2}^{{\wt}b-q}
Y_{M}(a,z_{0}+z_{2})Y_{M}(b,z_{2})
\end{eqnarray*}
on $M(U)(n)$, where $q=-1+l_{3}+\de_{i_{3}}(r)+r/T$.
\end{lem}

\pf Assume that $b\in V^s$. Take $v\otimes w\in
A_{g,n,m}(V)\otimes_{A_{g,m}(V)}U=M(U)(n).$ Then
\begin{eqnarray*} & &{\rm
Res}_{z_{0}}z_{0}^{-l}(z_{2}+z_{0})^{{\wt}a+q}z_{2}^{{\wt}b-q}Y_{M}(Y(a,z_{0})b,z_{2})(v\otimes w)\\
& &=\sum_{j\in{\mathbb Z}_{+}}{{\wt}a+q\choose
j}z_{2}^{{\wt}a+{\wt}b-j}Y_{M}(a_{j-l}b,z_{2})(v\otimes w)\\
& &=\sum_{j\in{\mathbb Z}_{+}}{{\wt}a+q\choose
j}\sum\limits_{k\in{\mathbb
Z}_{+}+\overline{i_{3}-r-s}/T}z_{2}^{-l+k-n+1}(a_{j-l}b)_{{\wt}a+{\wt}b-j+l-2-k+n}(v\otimes w)\\
& &=\sum_{k\in{\mathbb
Z}_{+}+\overline{i_{3}-r-s}/T}z_{2}^{-l+k-n+1}\sum_{j\in{\mathbb
Z}_{+}}{{\wt}a+q\choose j}\left((a_{j-l}b)\ast_{g,m,n}^{k}v\right)\otimes w\\
& &=\sum_{k\in{\mathbb
Z}_{+}+\overline{i_{3}-r-s}/T}z_{2}^{-l+k-n+1}\left(({\rm
Res}_{z}\frac{(1+z)^{{\wt}a+q}}{z^{l}}Y(a,z)b)\ast_{g,m,n}^{k}v\right)\otimes
w
\end{eqnarray*}
On the other hand, we have
\begin{eqnarray*}
& &{\rm
Res}_{z_{0}}z_{0}^{-l}(z_{0}+z_{2})^{{\wt}a+q}z_{2}^{{\wt}b-q}Y_{M}(a,z_{0}+z_{2})Y_{M}(b,z_{2})(v\otimes
w)\\
& &=\sum\limits_{i\in{\mathbb Z}_{+}}{-l\choose
i}(-1)^{i}a_{{\wt}a+q-l-i}z_{2}^{{\wt}b-q+i}Y_{M}(b,z_{2})(v\otimes w)\\
& &=\sum\limits_{i\in{\mathbb Z}_{+}}{-l\choose
i}(-1)^{i}a_{{\wt}a+q-l-i}\sum\limits_{j\geq
-n,j\in\Z-s/T}z_2^{-q+i+j}b_{{\wt}b-1-j}(v\otimes w)\\
& &=\sum_{i\in{\mathbb Z}_{+}}\sum_{\stackrel{j\geq
-n,j\in\Z-s/T}{l+i+j\geq 1+q-n}}{-l\choose i}(-1)^{i}z_2^{-q+i+j}
\left(a\ast_{g,m,j+n}^{l+i+j-1-q+n}
(b\ast_{g,m,n}^{j+n}v)\right)\otimes w\\
& &=\sum\limits_{k\in{\mathbb
Z}_{+}+\overline{i_{3}-r-s}/T}\sum\limits_{\stackrel{j\in\Z-s/T}{-n\leq
j\leq k+1+q-n-l}}z_2^{-l+k-n+1}(-1)^{k+1+q-n-j-l}{-l\choose
k+1+q-n-j-l}\\ & & \ \ \ \ \ \left((a\ast_{g,n,j+n}^{k}b)\ast_{g,m,n}^{k}v\right)\otimes w\\
& &=\sum\limits_{k\in{\mathbb
Z}_{+}+\overline{i_{3}-r-s}/T}\sum\limits_{\stackrel{j\in\Z+(i_{3}-s)/T}{0\leq
j\leq k+1+q-l}}z_2^{-l+k-n+1}(-1)^{k+1+q-j-l}{-l\choose k+1+q-j-l}\\
& & \ \ \ \ \
\left((a\ast_{g,n,j}^{k}b)\ast_{g,m,n}^{k}v\right)\otimes w.
\end{eqnarray*}
So it is enough to prove that
$$\sum\limits_{\stackrel{j\in\Z+(i_{3}-s)/T}{0\leq j\leq
k+1+q-l}}(-1)^{k+1+q-j-l}{-l\choose
k+1+q-j-l}a\ast_{g,n,j}^{k}b={\rm
Res}_{z}\frac{(1+z)^{{\wt}a+q}}{z^{l}}Y(a,z)b,
$$
for $k\geq 0$.
 Let $j=p+(i_{3}-s)/T-\de_{i_{3}}(s)+1$, $k=l_{1}+(i_{3}-r-s)/T$. Note that  $q=-1+l_{3}+\de_{i_{3}}(r)+r/T$,
 so
\begin{eqnarray*}
& &\sum\limits_{\stackrel{j\in\Z+(i_{3}-s)/T}{0\leq j\leq
k+1+q-l}}(-1)^{k+1+q-j-l}{-l\choose k+1+q-j-l}a\ast_{g,n,j}^{k}b\\
&
&=\sum\limits_{p=0}^{l_{1}+l_{3}+\de_{i_{3}}(r)+\de_{i_{3}}(s)-1-l}
(-1)^{l_{1}+l_{3}+\de_{i_{3}}(r)+\de_{i_{3}}(s)-1-l-p}{-l\choose
l_{1}+l_{3}+\de_{i_{3}}(r)+\de_{i_{3}}(s)-1-l-p}
\\& & \cdot\sum\limits_{i=0}^{p}(-1)^{i}{l_{1}+l_{3}+\de_{i_{3}}(r)+\de_{i_{3}}(s)-2-p+i\choose i} {\rm
Res}_{z}\frac{(1+z)^{{\wt}a+q}}{z^{l_{1}+l_{3}+\de_{i_{3}}(r)+\de_{i_{3}}(s)-1-p+i}}Y(a,z)b\\
&
&=\sum\limits_{p=0}^{l_{1}+l_{3}+\de_{i_{3}}(r)+\de_{i_{3}}(s)-1-l}(-1)^{p}{-l\choose
p}\sum\limits_{i=0}^{l_{1}+l_{3}+\de_{i_{3}}(r)+\de_{i_{3}}(s)-1-l-p}(-1)^{i}{l+p+i-1\choose
i}\\
& & \ \ \ \ \ \ {\rm
Res}_{z}\frac{(1+z)^{{\wt}a+q}}{z^{l+p+i}}Y(a,z)b.
\end{eqnarray*}
By Proposition 5.3 in \cite{DJ1}, we have
$$\sum\limits_{p=0}^{l_{1}+l_{3}+\de_{i_{3}}(r)+\de_{i_{3}}(s)-1-l}(-1)^{p}{-l\choose
p}\sum\limits_{i=0}^{l_{1}+l_{3}+\de_{i_{3}}(r)+\de_{i_{3}}(s)-1-l-p}(-1)^{i}{l+p+i-1\choose
i}\frac{1}{z^{l+p+i}} =\frac{1}{z^{l}}. $$
This finishes the proof. \qed

\begin{coro}\label{c05.10} Let $n=l_{3}+i_{3}/T\in\frac{1}{T}\Z_{+}$  with $l_{3}\in\Z_{+}$ and $0\leq i_{3}<T$.
Then for $a\in V^r$, we have
\begin{eqnarray*}
(z_{2}+z_{0})^{{\wt}a+q}Y_{M}(Y(a,z_{0})b,z_{2})=(z_{0}+z_{2})^{{\wt}a+q}Y_{M}(a,z_{0}+z_{2})Y_{M}(b,z_{2})
\end{eqnarray*}
on $M(U)(n),$ where $q=-1+l_{3}+\de_{i_{3}}(r)+r/T$.
\end{coro}

\begin{theorem}\label{t5.11} Let $U$ be an $A_{g,m}(V)$-module. Then
$M(U)=\bigoplus_{n\in\frac{1}{T}{\mathbb
Z}_{+}}A_{g,n,m}(V)\otimes_{A_{g,m}(V)} U$ is an admissible
$g$-twisted $V$-module with
$M(U)(n)=A_{g,n,m}(V)\otimes_{A_{g,m}(V)} U$ and with the
following universal property: for any weak $g$-twisted $V$-module
$W$ and any $A_{g,m}(V)$-morphism $\sigma: U\rightarrow
\Omega_{m}(W)$, there is a unique homomorphism $\bar\sigma:
M(U)\rightarrow W$ of weak $g$-twisted $V$-modules which extends
$\sigma$. Moreover, if $U$ cannot factor through $A_{g,m-1/T}(V)$
then $M(U)(0)\ne 0.$
\end{theorem}

It is clear from Theorem \ref{t5.11} that $M(U)$ is isomorphic to
the  $\bar{M}(U)$ given in Theorem \ref{t2.4}. We call $M(U)$ the
Verma type admissible $g$-twisted $V$-module generated by an
$A_{g,m}(V)$-module $U$.

\section{ $g$-rationality}
\def\theequation{6.\arabic{equation}}
\setcounter{equation}{0}

We use the bimodule theory developed in the previous sections
to prove another main theorem in this paper. That is, $V$ is $g$-rational
if and only if $A_{g}(V)$ is semisimple and any irreducible admissible
$g$-twisted $V$-module is ordinary. In the case $g=1,$ this result has
been obtained previously in \cite{DJ2} and \cite{DJ3}.

We need several lemmas. Let $A$ be an associative algebra and $U$
a left $A$-module. It is well known that the linear dual
$U^*=\Hom_{\C}(U,\C)$ is naturally a right $A$-module such that
$(fa)(u)=f(au)$ for $a\in A,$ $f\in U^*$ and $u\in U.$

\begin{lem}\label{l6.1} Let $V$ be a vertex operator algebra.
Assume that $A_{g}(V)$ is semisimple and $U^{i}$ for
$i=1,\cdots,s$ are all the inequivalent irreducible
$A_{g}(V)$-modules. Let
${M}(U^{i})=\bigoplus_{n\in\frac{1}{T}{\mathbb
Z}_{+}}{M}(U^{i})(n)$ be the Verma type admissible $g$-twisted
$V$-module generated by $A_{g}(V)$-module $U^{i}$. Then as an
$A_{g,n}(V)$-$A_{g}(V)$-bimodule,
$$
A_{g,n,0}(V)\cong
\bigoplus_{i=1}^{s}{M}(U^{i})(n)\otimes(U^{i})^{*}.$$
\end{lem}

\pf The proof is similar to that of Lemma 3.1 in \cite{DJ3}. \qed

\begin{lem}\label{l6.5}  Let $V$ be a  simple vertex operator algebra
such that $A_g(V)$ is finite dimensional. Then there exists
$N\in\frac{1}{T}\Z_{+}$ such that for any irreducible
$A_{g}(V)$-module $U$ and the irreducible admissible $g$-twisted
$V$-module $L(U)=\sum_{m\in\frac{1}{T}\Z_{+}}L(U)(m)$ generated by
$U$, $L(U)(n)\neq 0$, for all $n\in\frac{1}{T}\Z_{+}, n>N$.
\end{lem}

\pf  Since $A_g(V)$ is finite dimensional,
 there are only finitely many irreducible
admissible $g$-twisted $V$-modules. So it suffices to prove that
for the irreducible admissible $g$-twisted
$V$-module $W=\bigoplus_{m\in\frac{1}{T}\Z_{+}}W(m)$
there exists $N\in\frac{1}{T}\Z_{+}$ such that $W(n)\ne 0$
 for all $n\in\frac{1}{T}\Z_{+}, n>N$.

For $i=0,\cdots,T-1$ we set
$W^i=\bigoplus_{m\in\Z_+}W(m+\frac{i}{T}).$ Note from \cite{DM}
and \cite{DLM0} that  the $g$-invariants $V^{\langle g\rangle}$
again is a simple vertex operator algebra. From the definition of
admissible $g$-twisted $V$-modules we see that each $W^i$ is an
irreducible admissible $V^{\langle g\rangle}$-module (see
\cite{DY}). It is clear that for each $i$ there exists $n_i\geq 0$
such that $W^i(s)\ne 0$ if $s>n_i$ and $s\in \frac{i}{T}+\Z_+.$
Take $N$ to be the maximum of $n_i$ for $i=0,\cdots,T-1$ and the
lemma follows. \qed

For $m,n,p\in\frac{1}{T}{\mathbb Z}_{+}$, let
$$A_{g,n,p}(V)\ast_{g,m,p}^{n}A_{g,p,m}(V)=\{a\ast_{g,m,p}^{n}b|a\in
A_{g,n,p}(V),b\in A_{g,p,m}(V)\}.$$ Then
$A_{g,n,p}(V)\ast_{g,m,p}^{n}A_{g,p,m}(V)$ is an
$A_{g,n}(V)-A_{g,m}(V)$-subbimodule of $A_{g,n,m}(V)$
by Proposition \ref{p4.4}.

\begin{lem}\label{l6.6} Let $V$ be a  simple vertex operator algebra such that
$A_{g}(V)$ is
semisimple. Then there exists $N\in\frac{1}{T}\Z_{+}$ such that
$$A_{g,0,m}(V)\ast_{g,0,m}^{0}A_{g,m,0}(V)=A_{g}(V)$$
for all $m\in\frac{1}{T}\Z_{+}, m>N$.
\end{lem}

 \pf Let $N$ be the same as in Lemma \ref{l6.5}. For any $n\in\frac{1}{T}\Z_{+}$, it is easy to see
 that
$A_{g,0,n}(V)\ast_{g,0,n}^{0}A_{g,n,0}(V)$ is a two-sided ideal of
$A_{g}(V)$. Let $U$ be an irreducible module of $A_{g}(V)$ and
suppose that for some $m\in\frac{1}{T}\Z_{+}$, $m>N$,
\begin{equation}\label{ea6.0}A_{g,0,m}(V)\ast_{g,0,m}^{0}A_{g,m,0}(V)\otimes U=0.
\end{equation}
 Let $M(U)$ be
the Verma type admissible $g$-twisted $V$-module generated by $U$
and $M'(U)$ the maximal proper admissible $g$-twisted submodule of
$M(U)$.   Similar to the proof of Proposition 4.5.6 of \cite{LL}
(see also \cite{DM}), we have
$$
M''(U)=span\{u_{p}w|u\in V, p\in\Q, w\in M(U)(m)\}$$ is an
admissible $g$-twisted $V$-submodule of $M(U)$ generated by
$M(U)(m)$.
 By (\ref{ea6.0}) we know
that $M''(U)(0)=0$. So $M''(U)$ is a proper admissible $g$-twisted
$V$-submodule of $M(U)$ and $M''(U)\subset M'(U).$

Let $W(U)=M(U)/M'(U)$, then $W(U)$ is the irreducible admissible
$g$-twisted  $V$-module generated by $U$ and $W(U)(m)=0$. This is
in contradiction with Lemma \ref{l6.5}. Thus for all
$m\in\frac{1}{T}\Z_{+}$, $m>N$,
$$A_{g,0,m}(V)\ast_{g,0,m}^{0}A_{g,m,0}(V)\otimes U\cong U.$$
Now the lemma follows from Lemma \ref{l6.1}. \qed

Recall from  Proposition \ref{p4.4} that $\varphi$ :
$A_{g,n,p}(V)\otimes_{A_{g,p}(V)}A_{g,p,m}(V)\rightarrow
A_{g,n,m}(V)$ is an $A_{g,n}(V)-A_{g,m}(V)$-bimodule homomorphism
defined by
 $$
\varphi(u\otimes v)=u\ast_{g,m,p}^{n}v,$$ for $u\in A_{g,n,p}(V)$,
$v\in A_{g,p,m}(V)$ and $m,p,n\in\frac{1}{T}\Z_{+}.$ The following
lemma is an immediate consequence of Lemma \ref{l6.6}.
\begin{lem}\label{l6.7} Let $V$
be a  simple vertex operator algebra such that  $A_{g}(V)$ is
semisimple. Then there exists $N\in\frac{1}{T}\Z_{+}$ such that
the $A_{g}(V)-A_{g}(V)$-bimodule homomorphism $\varphi$ from
$A_{g,0,n}(V)\otimes_{A_{g,n}(V)}A_{g,n,0}(V)$ to $A_{g}(V)$ is an
isomorphism for each $n\in\frac{1}{T}\Z_{+}$, $n>N$.
\end{lem}

\begin{theorem}\label{t6.8} Let $V$
be a simple  vertex operator algebra such that $A_{g}(V)$ is
semisimple. Let $U$ be an irreducible module of $A_{g}(V)$, then
the Verma type admissible $g$-twisted $V$-module
$M(U)=\bigoplus_{n\in\frac{1}{T}{\mathbb
Z}_{+}}A_{g,n,0}(V)\otimes_{A_{g}(V)}U$ generated by $U$ is
irreducible.
\end{theorem}

\pf The same proof of Theorem 3.4 of \cite{DJ3} is valid here. \qed

We have already mentioned that the Verma type admissible
$g$-twisted $V$-module $M(U)$ generated by an irreducible
$A_{g}(V)$-module $U$ in general is not irreducible. The
assumption that $A_g(V)$ is semisimple is crucial. This result is
a foundation of Theorem \ref{t6.9} below.

As in \cite{DJ3} we now introduce an invariant bilinear
pairing $(\cdot,\cdot)$ on $M(U^{*})\times M(U)$, for an
$A_{g,m}(V)$-module $U$, which is an analogue of the contravariant
forms for Verma modules over 
Kac-Moody Lie algebras or the Virasoro algebra. This
bilinear pairing will also be helpful to the proof of Theorem
\ref{t6.9}.

Let $m\in\frac{1}{T}\Z_{+}$. Recall from Theorem \ref{t2.4} (8)
and Proposition \ref{p4.2} that the linear map $\phi$:
$V\rightarrow V$ defined by $\phi(v)=e^{L(1)}(-1)^{L(0)}v$ for
$v\in V$ induces an anti-isomorphism from $A_{g^{-1},m}(V)$ to
$A_{g,m}(V)$ and a linear isomorphism from $A_{g^{-1},n,m}(V)$ to
$A_{g,m,n}(V)$ such that
$\phi(a*_{g^{-1},m,p}^{n}b)=\phi(b)*_{g,n,p}^{m}\phi(a)$, for
$a\in A_{g^{-1},n,p}(V)$, $b\in A_{g^{-1},p,m}(V)$,
$p\in\frac{1}{T}\Z_{+}$.

 Let $U$ be an $A_{g,m}(V)$-module and $U^{*}$ the
dual space of $U$. Then $U^*$ is an $A_{g^{-1},m}(V)$-module such that
$$
(u\cdot f)(x)=f(\phi(u)\cdot x)=(f, \phi(u)\cdot x)$$ for $u\in
A_{g^{-1},m}(V)$, $f\in U^*$ and $x\in U$.

 Now let $M(U)=\bigoplus_{n\in\frac{1}{T}{\mathbb
Z}_{+}}A_{g,n,m}(V)\otimes_{A_{g,m}(V)} U$ and
$M(U^*)=\bigoplus_{n\in\frac{1}{T}{\mathbb
Z}_{+}}A_{g^{-1},n,m}(V)$ $\otimes_{A_{g^{-1},m}(V)} U^*$ be the
Verma type admissible $g$-twisted and $g^{-1}$-twisted $V$-modules
generated by $U$ and $U^*$ respectively. We define a bilinear
pairing $(\cdot,\cdot)$ on $M(U^*)\times M(U)$ as follows:
$$
(x\otimes f, y\otimes u)=(f, [(\phi(x)\ast_{g,m,n}^{m} y)]\cdot
u)$$ 
for $x\in A_{g^{-1},n,m}(V)$, $y\in A_{g,n,m}(V), f\in U^*,
u\in U,n\in\frac{1}{T}{\mathbb Z}_{+}$; and
$$
(A_{g^{-1},p,m}(V)\otimes_{A_{g^{-1},m}(V)}U^*,A_{g,n,m}(V)\otimes_{A_{g,m}(V)}U)=0$$
 for $p\neq n$. That is $(M(U^*)(p),M(U)(n))=0$ if $p\ne n.$

As in \cite{DJ3} we have
\begin{prop} \label{p6.10} Let $U$ be an $A_{g,m}(V)$-module 
for $m\in\frac{1}{T}\Z_+.$ Then

(1) The bilinear pairing
$(\cdot,\cdot)$ on $M(U^*)\times M(U)$ is well defined 
and is invariant in the sense
that
$$
(Y_{M(U^*)}(u,z)w',w)=(w',
Y_{M(U)}(e^{zL(1)}(-z^{-2})^{L(0)}u,z^{-1})w)
$$
for $w'\in M(U^*) ,w\in M(U),$ and $u\in V$.

(2) The space 
$$J(U)=\{w\in M(U)|(w', w)=0, w'\in M(U^*)\}
$$
is the maximal proper admissible $g$-twisted $V$-submodule of
$M(U)$ such that $$J(U)\cap M(U)(m)=0.$$ In particular, if $U$ is
irreducible then $J(U)$ is the unique maximal proper admissible
$g$-twisted  submodule of $M(U).$

(3) Let $V$ be a  simple vertex operator algebra  such that
$A_{g}(V)$ is semisimple. Let $U$ be an irreducible
$A_{g}(V)$-module, then the bilinear pairing $(\cdot,\cdot)$ on
$M(U^*)\times M(U)$ is non-degenerate.
\end{prop}

The analogue of Theorem 5.3 of \cite{DJ3} is the following -- the
second main theorem in this paper with the similar proof.

\begin{theorem}\label{t6.9} Let $V$ be a simple vertex operator
algebra and $g$ an  automorphism of $V$ of finite order. Then $V$
is $g$-rational if and only if $A_{g}(V)$ is semisimple and each
irreducible admissible $g$-twisted $V$-module is ordinary.
\end{theorem}

We remark that the condition that each irreducible admissible
$g$-twisted $V$-module is ordinary
 holds for all known simple vertex
operator algebras  and finite order automorphisms. Although we
firmly believe that this is true in general, we can not prove this
in this paper.

\end{document}